\documentclass[11pt,a4paper]{article}
\usepackage[latin1]{inputenc}
\usepackage[english]{babel}
\usepackage{amsmath,amsthm}
\usepackage{amsfonts}
\usepackage{mathrsfs} 
\usepackage{latexsym}
\usepackage[dvips]{graphicx}
\usepackage{float}
\usepackage[natural]{xcolor}
\usepackage{algorithm}
\usepackage{algorithmicx}
\usepackage{algpseudocode}
\usepackage{slashed}
\usepackage{enumerate,enumitem}
\usepackage{multirow}

\usepackage[colorlinks,linkcolor=blue]{hyperref}
\usepackage[pagewise]{lineno}
\usepackage{fullpage}
\usepackage[title]{appendix}
\usepackage{todonotes}
\usepackage{amssymb}
\usepackage{slashed}
\usepackage{tikz}
\tikzstyle{aNode} = [circle, fill = black]
\tikzstyle{bNode} = [circle,draw = black, thick]
\usepackage{subcaption}
\usepackage{stmaryrd}
\usepackage{listings}
\usepackage{bm}
\usepackage{bbm}
\newcounter{propcounter}
\usepackage{enumitem}

\newcommand{\sub}{\subseteq}
\newcommand{\eps}{\varepsilon}

\makeatletter
\renewcommand{\fnum@figure}{\small\textbf{\figurename~\thefigure}}

\makeatother

\theoremstyle{plain}
\newtheorem{thm}{Theorem}[section]
\newtheorem{theorem}[thm]{Theorem}
\newtheorem{conjecture}[thm]{Conjecture}
\newtheorem{lemma}[thm]{Lemma}
\newtheorem{corollary}[thm]{Corollary}
\newtheorem{proposition}[thm]{Proposition}

\theoremstyle{definition}

\newtheorem{question}[thm]{Question}
\newtheorem{problem}[thm]{Problem}
\newtheorem{remark}{Remark}

\newtheorem{definition}[thm]{Definition}
\newtheorem{claim}[thm]{Claim}
\newtheorem{fact}[thm]{Fact}

\newtheorem{example}[thm]{Example}

\newtheorem{observation}[thm]{Observation}
\newtheorem{defn-thm}[thm]{Definition-Theorem}
\numberwithin{equation}{section}

\newcommand{\btheorem}{\begin{theorem}}
	\newcommand{\etheorem}{\end{theorem}}
\newcommand{\bconjecture}{\begin{conjecture}}
	\newcommand{\econjecture}{\end{conjecture}}
\newcommand{\bproposition}{\begin{proposition}}
	\newcommand{\eproposition}{\end{proposition}}
\newcommand{\bdefinition}{\begin{definition}}
	\newcommand{\edefinition}{\end{definition}}
\newcommand{\bcorollary}{\begin{corollary}}
	\newcommand{\ecorollary}{\end{corollary}}
\newcommand{\pr}{\begin{proof}}
	\newcommand{\oof}{\end{proof}}
\newcommand{\bclaim}{\begin{claim}}
	\newcommand{\eclaim}{\end{claim}}
\newcommand{\bquestion}{\begin{question}}
	\newcommand{\equestion}{\end{question}}
\newcommand{\bfact}{\begin{fact}}
	\newcommand{\efact}{\end{fact}}
\newcommand{\bremark}{\begin{remark}}
	\newcommand{\eremark}{\end{remark}}
\newcommand{\eexample}{\end{example}}
\newcommand{\bexample}{\begin{example}}
\newcommand{\ma}{\end{lemma}}
\newcommand{\lem}{\begin{lemma}}

\newcommand{\ppoints}[1]{%
	\begin{tikzpicture}[inner sep = 0.7pt, #1]%
		\node (1) at (0,-2) [aNode]{};
		\node (3) at (1.5,-2) [aNode]{};
		\node (2) at (0.75,-1) [aNode]{};
	\end{tikzpicture}%
}

\def\points{\ppoints{scale=0.11}}

\begin{document}

\title{Hypergraphs with a quarter uniform Tur\'an density}
\author{Hao Li \thanks{Laboratoire Interdisciplinaire des Sciences du Num\'{e}rique, UMR9015 CNRS-Universit\'{e} Paris-Saclay, 1 Rue Raimond Castaing, 91190 Gif-sur-Yvette, France. Email: \texttt{li@lri.fr}.}
	\quad Hao Lin \thanks{Faculty of Electrical Engineering, Mathematics and Computer Science, Delft University of Technology, Delft, Netherlands. Email: \texttt{haolinz6@qq.com}.}
	\quad Guanghui Wang\thanks{School of Mathematics, Shandong University, Jinan, China. Email: \texttt{ghwang@sdu.edu.cn}.}
	\quad Wenling Zhou\textsuperscript{\P}\thanks{School of Mathematics, Shandong University, Jinan, China. Email: \texttt{gracezhou@sdu.edu.cn}. \\
		\text{\textsuperscript{\P}Corresponding author.}}}
	
\date{}
\maketitle

\begin{abstract}
The \emph{uniform Tur\'an density $\pi_{\points}(F)$} of a ($3$-uniform) hypergraph $F$ is the supremum of $d$ for which there are infinitely many $F$-free hypergraphs with the property that 
every induced subhypergraph of $H$ on a  linearly sized vertex set has edge density at least $d$. Determining $\pi_{\points}(F)$ for given hypergraphs $F$ was suggested by Erd\H{o}s and S\'os in the 1980s. However, there are very few hypergraphs whose uniform Tur\'an density has been determined.

In this paper, we are the first to establish a verifiable condition for hypergraphs $F$ with $\pi_{\points}(F) = 1/4$.
 In particular, currently known hypergraphs  whose uniform Tur\'an density is $1/4$, such as $K_4^{(3)-}$ studied in [Israel J. Math. 211 (2016), 349-366] and [J. Eur. Math. Soc. 20 (2018), 1139-1159], and $F^{\star}_5$ studied in  [arXiv:2211.12747], satisfy this condition. Moreover, we also identify some new  hypergraphs whose uniform Tur\'an density is also $1/4$.
 \smallskip
 
 \noindent{\bf Keywords:} Uniform Tur\'an density, Quasi-random hypergraphs,  Hypergraph regularity method
\end{abstract}

\section{Introduction}
	For a positive integer $n$, we denote by $[n]$ the set $\{1,\dots,n\}$.
	Given $k\geq 2$, for a finite set $V$, we use $[V]^k$ to denote the collection of all subsets of $V$ of size $k$.
	We may drop one pair of brackets and write $[n]^k$ instead of $[[n]]^k$. In addition,  we simply write $i_1\dots i_k$ to denote a set  $\{i_1, \dots, i_k\} \subset \mathbb Z$ with $i_1<\dots<i_k$.
	A {\it $k$-uniform hypergraph} $H$ (or \textit{$k$-graph} for short)  is a pair $H=(V,E)$, where $V$ is a finite set of \textit{vertices} and $E\subseteq [V]^{k}$ is a set of \textit{edges} (or \emph{$k$-edges}).
	Given a $3$-graph $F$ and a vertex $v\in V(F)$, let
	$\partial F:=\{\{x, y\}\subseteq [V(F)]^2: \text{there is}~ z\in V(F)  \text{~such that~}  \{x,y, z\}\in E(F)\}$
	be the \emph{shadow} of $F$, and $L_v(F):=\{\{x, y\}\subseteq [V(F)]^2:\{x,y,v\}\in E(F)\}$ be the link graph of $v$. Clearly, $\partial F=\cup_{v\in V(F)}L_v(F)$. Given two $3$-graphs $F$ and $H$, we write $F\subseteq H$ to denote that $H$ contains a copy of $F$ as a subhypergraph.

	A classical extremal problem for (hyper)graphs introduced by Tur\'an~\cite{Tu-graph} about 80 years ago,  
	asks to study for a given (hyper)graph $F$ its \textit{Tur\'an number} ${\rm ex}(n, F)$, the maximum number of edges in an $F$-free (hyper)graph on $n$ vertices.
	It is a long-standing open problem in Extremal Combinatorics to develop some understanding of these numbers. Ideally, one would like to compute them exactly, but even asymptotic results (for $k\ge 3$) are currently only known in certain cases. For more details, see the excellent survey by Keevash~\cite{Keevash-survey}.
	It is well known and not hard to observe that the sequence ${\rm ex}(n, F)/\binom{n}{k}$ is decreasing.
	Thus, 
	one often focuses on the \textit{Tur\'an density} $\pi(F)$ of $F$ defined by
	\[
	\pi(F)=\lim_{n\to \infty}\frac{{\rm ex}(n,F)}{\binom{n}{k}}.
	\] 
	Tur\'an densities are well understood for graphs, i.e., $2$-graphs.
	Indeed, the Mantel's theorem~\cite{mantel1907problem} and the  Tur\'an's theorem~\cite{Tu-graph} give  the Tur\'an numbers of complete graphs exactly,
	and Erd\H{o}s and Stone~\cite{E-Stone} (also see Erd\H{o}s and Simonovits~\cite{E-Simonovits}) determined the Tur\'an density of any graph $F$ to be equal to 
	$\frac{\chi(F)-2}{\chi(F)-1}$,
	where $\chi(F)$ denotes the {\it chromatic number} of $F$,
	that is the minimum number of colors used to color $V(F)$ such that any two adjacent vertices receive distinct colors. 
	In particular, let $\Pi^{(k)}=\{\pi(F): F \text{~is~a~} k\text{-graph}\}$. Then
	\[
	\Pi^{(2)}=\{0, \frac{1}{2}, \frac{2}{3}, \dots, \frac{t-2}{t-1},\dots\}.
	\]
	However, the extremal problems for hypergraphs are notoriously difficult, even for ``simple" hypergraphs like the $3$-uniform clique $K^{(3)}_4$ on four vertices and  $K^{(3)-}_4$ the $3$-graph with four vertices and three edges. Despite much effort and attempts so far, our knowledge is somewhat limited, for example
	the best-known  bounds (see~\cite{k43lowbound, k43-lowbound, k43upbound, k43-upbound}) for $\pi(K^{(3)}_4)$ and $\pi(K^{(3-)}_4)$ are
	\[
	0.\dot 5=\frac{5}{9}\le \pi(K^{(3)}_4)\le 0.5615 ~~~~\text{and}~~~~ 0.\dot 2 \dot 8 \dot 5 \dot 7 \dot 1 \dot 4=\frac{2}{7}\le \pi(K^{(3)-}_4)\le 0.2871.
	\]
	Moreover, it is known that the corresponding set  $\Pi^{(3)}$ is much more complicated than $\Pi^{(2)}$ and as a subset of the reals it is not well-ordered (see, e.g.,~\cite{frankl1984hypergraphs, Pikhurko2014, tight-cycles-2022tur}).

	\subsection{Tur\'an problems in  uniformly dense hypergraphs}
	Most known and conjectured extremal constructions for Tur\'an problems contain large independent sets, i.e., linear-size sets of vertices without edges. 
	In 1982, Erd\H{o}s and S\'os~\cite{E-Sos} suggested a variation restricting Tur\'an problems only to those $F$-free $3$-graphs that are uniformly dense on large subsets of the vertices, formally defined as follows. 
	Given real numbers $d \in [0, 1]$ and $\mu>0$, a 3-graph $H=(V, E)$ is said to be {\it $(d, \mu)$-dense} (or~{\it uniformly dense}) if for any $U\subseteq V$, it holds that
	\begin{equation}\label{eq:vertex-dense}
		\left|[U]^3\cap E\right|\ge d\binom{|U|}{3}-\mu|V|^3.
	\end{equation}
	Specifically, restricting to uniformly dense hypergraphs, the \emph{uniform Tur\'an density} $\pi_{\points}(F)$ for a given
	$3$-graph $F$ is defined as
	\begin{equation*}
		\begin{split}
			\pi_{\points}(F) = \sup \{ d\in [0,1] &: \text{for\ every\ } \mu>0 \ \text{and\ } n_0\in \mathbb{N},\ \text{there\ exists\ an\ } F \text{-free} \\
			&\quad (d,\mu)\text{-dense}~  \text{$3$-graph~} H\ \text{with~} |V(H)|\geq n_0 \}.
		\end{split}
	\end{equation*}
	With this notation at hand, Erd\H{o}s and S\'os  asked to determine $\pi_{\points}(K^{(3)-}_4)$ and $\pi_{\points}(K^{(3)}_4)$. 
	However, determining $\pi_{\points}(F)$  of a given $3$-graph $F$ is also very challenging. $\pi_{\points}(K^{(3)-}_4)=1/4$ was solved recently by Glebov, Kr\'al' and Volec~\cite{k43minus-1}, and independently by Reiher, R\"odl and
	Schacht~\cite{k43minus-2}.
	The  conjecture for $\pi_{\points}(K^{(3)}_4)=1/2$  has been an urgent problem in this area since R\"odl~\cite{k43-rodl} gave a quasi-random construction in 1986. 
	In addition to $K^{(3)-}_4$, Reiher, R\"odl and
	Schacht~\cite{vanishing} characterised
	all $3$-graphs $F$ with $\pi_{\points}(F)=0$.

	\begin{theorem}[{\cite[Theorem 1.2]{vanishing}}]\label{0-density}
		A $3$-graph $F$ on $f$ vertices satisfies $\pi_{\points}(F)=0$ if and only if it has the following property.
		There exists an  ordering $(v_1, \dots, v_f)$ of $V(F)$ and there is a  $3$-coloring $\chi: \partial F\to \{ {\color {red}red}, {\color{blue}blue}, {\color{green}{green}} \}$  
		such that every edge $\{v_i, v_j, v_k\}$ of $F$ with $i<j< k$ satisfies
		\begin{equation}\label{eq-0-density}
			\left(\chi(v_i,v_j), \chi(v_j,v_k), \chi(v_i, v_k)\right)=({\color {red}red},  {\color{blue}blue},  {\color{green}green}).
		\end{equation} 
	\end{theorem}
	
	Meanwhile, they~\cite{RRS-Mantel} also proposed the following problem.

	\begin{problem}[{\cite[Problem 1.7]{RRS-Mantel}}]
		Determine $\pi_{\points}(F)$ for all $3$-graphs $F$.
	\end{problem}

	Very recently, Buci{\'c}, Cooper, Kr\'al', Mohr and Munh\'a Correia~\cite{cycle-uniform} determined the uniform Tur\'an density of all tight $3$-uniform cycles of
	length at least five; Garbe, Kr\'al' and Lamaison~\cite{1/27} showed a specific family of $3$-graphs $F$ with $\pi_{\points}(F)=1/27$; Chen and Sch{\"u}lke~\cite{chen2022beyond} proved $\pi_{\points}(F^{\star}_5)=1/4$, where $F^{\star}_5$ is a $3$-graph on 5 vertices
	that is obtained from $K^{(3)-}_4$
	by adding a new vertex whose link graph forms a matching of size $2$
	on the vertices of $K^{(3)-}_4$. Given a $3$-graph $F$,  the \textit{$t$-blow-up} of $F$, denoted by $F(t)$, is the $3$-graph obtained from $F$ by replacing every vertex of $F$ by $t$ copies of itself. Furthermore, we say $F'$ is a \textit{blow-up} of $F$
	if $F\subseteq F'\subseteq F(t)$ for some positive integer $t$. By the Hypergraph Blow-up Lemma of Keevash~\cite{hypergraph-blow-up}, one can easily verify that  $\pi_{\points}(F)=\pi_{\points}(F')$ for every $3$-graph $F$ and its blow-up $F'$.
Aside from the aforementioned $3$-graphs and their blow-ups, there are currently very few results for the uniform Tur{\'a}n density of $3$-graphs (see Reiher's survey~\cite{reiher2020extremal}).

	\subsection{Our result}
	In this work, we first to establish a verifiable condition for $3$-graphs with the uniform Tur\'an  density $1/4$. We start with some notation.
	For simplicity, we define $\mathcal F^{\ge 1/4}$ as the collection of $3$-graphs $F$ with $f$ vertices that do not satisfy the following property:
	\begin{enumerate} 
		\item[{\rm ($\clubsuit$)}] There exists an  ordering $(v_1, \dots, v_{f})$ of $V(F)$  and there is a  $2$-coloring $\chi: \partial F\to \{{\color {red}red}, {\color{blue}blue}\}$  such that every edge $\{v_i, v_j, v_k\}$ of $F$ with $i<j< k$ satisfies
		\[
		\left(\chi(v_i,v_j), \chi(v_j,v_k), \chi(v_i,v_k)\right)=({\color {red}red}, {\color {red}red}, {\color{blue}blue}) \text{~or~} ({\color{blue}blue},  {\color{blue}blue}, {\color {red}red}).
		\]
	\end{enumerate}

	Using a probabilistic construction, we have the following fact.
	\begin{fact}\label{low-bound}
		For any $F\in \mathcal F^{\ge 1/4}$, we have $\pi_{\points}(F)\ge 1/4$.
	\end{fact}

	\begin{proof}
		Given a positive integer $n$, we consider a random $2$-coloring $\chi: [n]^2\to \{{\color {red}red}, {\color{blue}blue}\}$ with each color associated to a pair with probability $1/2$ independently.
		We define a random $3$-graph $H(n)$ on set $[n]$ as follows. For any triple $ijk\in [n]^3$,  $ijk \in E(H(n))$ if and only if $ijk$ satisfies 
		\[
		\left(\chi(i,j), \chi(j,k), \chi(i,k)\right)=({\color {red}red}, {\color {red}red}, {\color{blue}blue}) \text{~or~} ({\color{blue}blue},  {\color{blue}blue}, {\color {red}red}).
		\]
		Observe that for any $ijk\in [n]^3$, the probability of the event ``$ijk \in E(H(n))$" happening is $1/4$.
		By standard probabilistic arguments, we can easily check that for any fixed $\mu>0$,
		$H(n)$ is $(1/4, \mu)$-dense with probability tending to $1$ as $n\to \infty$.
		In addition, each subhypergraph of $H(n)$ satisfies the property {\rm ($\clubsuit$)}. Therefore,  $H(n)$ is $F$-free and $\pi_{\points}(F)\ge 1/4$ for all $F\in \mathcal F^{\ge 1/4}$.
	\end{proof}

	Next, we define $\mathcal F^{\le 1/4}$ as the collection of $3$-graphs $F$ with $f$ vertices satisfying the following property:

	\begin{itemize}
		\item[{\rm ($\spadesuit$)}]
		There exists an order $(v_1, \dots, v_{f})$ of $V(F)$ with an integer $i^\star\in [f]$, and
		there is a  $6$-coloring 
		$\chi: \partial F\to \{\textcolor{red}{red},\textcolor{blue}{blue}, \textcolor{green}{green},\textcolor{violet}{violet},\textcolor{cyan}{cyan}, \textcolor{black}{black}\}$
		such that $E(F)$ can be partitioned into three subsets $E_1$, $E_2$ and $E_3$  with the following properties:
		\stepcounter{propcounter}
		\begin{enumerate}[label = \rm({\bfseries \Alph{propcounter}\arabic{enumi}})]
			\item \label{p11} For $\{v_i,v_j,v_k\}\in E_1$, $\left(\chi(v_i,v_j), \chi(v_j,v_k), \chi(v_i,v_k)\right)=(\textcolor{blue}{blue}, \textcolor{violet}{violet}, \textcolor{red}{red})$ and $ijk\in [f]^3$;
			\item \label{p22} For $\{v_i,v_j,v_k\}\in E_2$, $\left(\chi(v_i,v_j), \chi(v_j,v_k),  \chi(v_i,v_k)\right)=(\textcolor{green}{green},  \textcolor{cyan}{cyan},\textcolor{blue}{blue})$ and $ijk\in [i^\star]^3$;
			\item \label{p33} For $\{v_i,v_j,v_k\}\in E_3$, $\left(\chi(v_i,v_j), \chi(v_j,v_k), \chi(v_i, v_k)\right)= (\textcolor{green}{green}, \textcolor{black}{black},\textcolor{red}{red})$ and $1\le i<j<i^\star< k\le f$.
		\end{enumerate}
	\end{itemize}
	
	Combining the properties {\rm ($\clubsuit$)} and  {\rm ($\spadesuit$)}, we obtain the following general result.
	
	\begin{theorem}\label{main-thm1}
		For any $F\in \mathcal F^{\ge 1/4}\cap\mathcal F^{\le 1/4}$, we have $\pi_{\points}(F)= 1/4$.
		
	\end{theorem}
	
	Through further analysis of the property {\rm ($\spadesuit$)}, we get the following further result.

\begin{theorem}\label{less-than1/4}
	Let $d^*\in \mathbb R $ denote the solution {\rm(}unique real root{\rm)} of the equation $(2-x)^3=27x$. For every $F \in \mathcal F^{\le 1/4}$, if an edge subset $E_i$ is empty for some $i\in [3]$, then $\pi_{\points}(F) \le d^*$.
	\end{theorem} 
	
	\begin{remark}
Note that $d^* <1/4$. 
Thus, Theorem~\ref{less-than1/4} implies that  for any $F\in  \mathcal F^{\le 1/4}$, if $\pi_{\points}(F)=1/4$, then  $E_1$, $E_2$ and $E_3$ must be non-empty in the property 
		{\rm ($\spadesuit$)}.
	\end{remark}

	\begin{figure}
		\centering
		\begin{subfigure}[b]{0.2\textwidth}
			\centering
			\begin{tikzpicture}
				\node[draw, circle, minimum size=0.5cm] (a) at (0.2,2.3)  {a};
				\node[draw, circle, minimum size=0.5cm] (b) at (0.2,0.3)  {b};
				\node[draw, circle, minimum size=0.5cm] (c) at (2.2,0.3)  {c};
				\node[draw, circle, minimum size=0.5cm] (d) at (2.2,2.3)  {d};
				\draw[line width=3pt,red,opacity=0.5,rounded corners=10pt](a) -- (b) -- (c) -- (a);
				\draw[line width=3pt,blue,opacity=0.5,rounded corners=10pt](0.3,2) -- (0.3,0.65);
				\draw[line width=3pt,blue,opacity=0.5,rounded corners=10pt] (0.4,0.6)--(d);
				\draw[line width=3pt,blue,opacity=0.5,rounded corners=10pt](d)--(a);
				\draw[line width=3pt,green,opacity=0.8,rounded corners=10pt](0.5,2.1) -- (2.05,0.55);
				\draw[line width=3pt,green,opacity=0.8,rounded corners=10pt](2.2,0.6)--(2.2,1.95);
				\draw[line width=3pt,green,opacity=0.8](1.85,2.2) --(0.5,2.2);
				\node[draw, circle] (x) at (1.1,3.5)  {x} ;
				\draw[line width=3pt,black,opacity=0.8,rounded corners=10pt](b)--(x) -- (a);
				\draw[line width=3pt,black,opacity=0.8,rounded corners=10pt](0.1,2) -- (0.1,0.66);
				\draw[line width=3pt,black,opacity=0.8,rounded corners=10pt](c)--(x) -- (d);
				\draw[line width=3pt,black,opacity=0.8,rounded corners=10pt](2.3,0.6)--(2.3,1.95);
			\end{tikzpicture}
			\subcaption{The $3$-graph $F^{\star}_5$.}\label{fig:picture1}
		\end{subfigure}
		\hfill
		\begin{subfigure}[b]{0.25\textwidth}
			\centering
			\begin{tikzpicture}
				\node[draw, circle] (a) at (0.2,2.3)  {a} ;
				\node[draw, circle] (b) at (0.2,0.3)  {b} ;
				\node[draw, circle] (c) at (2.2,0.3)  {c} ;
				\node[draw, circle] (d) at (2.2,2.3)  {d} ;
				\draw[line width=3pt,red,opacity=0.5,rounded corners=10pt](a) -- (b)--(c)--(a);
				\draw[line width=3pt,blue,opacity=0.5,rounded corners=10pt](0.3,2) -- (0.3,0.65);
				\draw[line width=3pt,blue,opacity=0.5,rounded corners=10pt] (0.4,0.6)--(d);
				\draw[line width=3pt,blue,opacity=0.5,rounded corners=10pt](d)--(a);
				\draw[line width=3pt,green,opacity=0.8,rounded corners=10pt](0.5,2.1) -- (2.05,0.55);
				\draw[line width=3pt,green,opacity=0.8,rounded corners=10pt](2.2,0.6)--(2.2,1.95);
				\draw[line width=3pt,green,opacity=0.8](1.85,2.2) --(0.5,2.2);
				\node[draw, circle] (x) at (1.1,3.5)  {x} ;
				\draw[line width=3pt,black,opacity=0.8,rounded corners=10pt](b)--(x) -- (a);
				\draw[line width=3pt,black,opacity=0.8,rounded corners=10pt](0.1,2) -- (0.1,0.66);
				\draw[line width=3pt,black,opacity=0.8,rounded corners=10pt](c)--(x) -- (d);
				\draw[line width=3pt,black,opacity=0.8,rounded corners=10pt](2.3,0.6)--(2.3,1.95);
				\node[draw, circle] (y) at (-1.2,1.2)  {y} ;
				\draw[line width=3pt,violet,opacity=0.5,rounded corners=10pt](c)--(y)--(a);
				\draw[line width=3pt,violet,opacity=0.5,rounded corners=10pt](0.35,2)-- (1.9,0.45);
				\draw[line width=3pt,violet,opacity=0.5,rounded corners=10pt](b)--(y) -- (d);
				\draw[line width=3pt,violet,opacity=0.5,rounded corners=10pt](0.3,0.65)--(1.9,2.15);
			\end{tikzpicture}
			\caption{The $3$-graph $F^{\star}_6$.}\label{fig:picture2}
		\end{subfigure}
		\hfill
		\begin{subfigure}[b]{0.4\textwidth}
			\centering
			\begin{tikzpicture}
				\node[draw, circle] (a) at (0.2,2.3)  {a} ;
				\node[draw, circle] (b) at (0.2,0.3)  {b} ;
				\node[draw, circle] (c) at (2.2,0.3)  {c} ;
				\node[draw, circle] (d) at (2.2,2.3)  {d} ;
				\draw[line width=3pt,red,opacity=0.5,rounded corners=10pt](a) -- (b)--(c)--(a);
				\draw[line width=3pt,blue,opacity=0.5,rounded corners=10pt](0.3,2) -- (0.3,0.65);
				\draw[line width=3pt,blue,opacity=0.5,rounded corners=10pt] (0.4,0.6)--(d);
				\draw[line width=3pt,blue,opacity=0.5,rounded corners=10pt](d)--(a);
				\draw[line width=3pt,green,opacity=0.8,rounded corners=10pt](0.5,2.1) -- (2.05,0.55);
				\draw[line width=3pt,green,opacity=0.8,rounded corners=10pt](2.2,0.6)--(2.2,1.95);
				\draw[line width=3pt,green,opacity=0.8](1.85,2.2) --(0.5,2.2);
				\node[draw, circle] (x) at (1.1,3.5)  {x} ;
				\draw[line width=3pt,black,opacity=0.8,rounded corners=10pt](b)--(x) -- (a);
				\draw[line width=3pt,black,opacity=0.8,rounded corners=10pt](0.1,2) -- (0.1,0.66);
				\draw[line width=3pt,black,opacity=0.8,rounded corners=10pt](c)--(x) -- (d);
				\draw[line width=3pt,black,opacity=0.8,rounded corners=10pt](2.3,0.6)--(2.3,1.95);
				\node[draw, circle] (y) at (-1.2,1.2)  {y} ;
				\draw[line width=3pt,violet,opacity=0.5,rounded corners=10pt](c)--(y)--(a);
				\draw[line width=3pt,violet,opacity=0.5,rounded corners=10pt](0.35,2)-- (1.9,0.45);
				\draw[line width=3pt,violet,opacity=0.5,rounded corners=10pt](b)--(y) -- (d);
				\draw[line width=3pt,violet,opacity=0.5,rounded corners=10pt](0.3,0.65)--(1.9,2.15);
				\node[draw, circle] (z) at (4,1.2)  {z} ;
				\draw[line width=3pt,yellow,opacity=0.8,rounded corners=10pt](a)--(z)--(d);
				\draw[line
				width=3pt,yellow,opacity=0.8,rounded corners=10pt](0.5,2.4)-- (1.85,2.4);
				\draw[line width=3pt,yellow,opacity=0.8,rounded corners=10pt](b)--(z) -- (c);
				\draw[line width=3pt,yellow,opacity=0.8,rounded corners=10pt](0.6,0.2)--(1.9,0.2);
			\end{tikzpicture}
			\caption{The $3$-graph $F^{\star}_7$.}\label{fig:picture3}
		\end{subfigure}
		\caption{}\label{fig:3-pictures}
	\end{figure}

We notice that the family $\mathcal F^{\ge 1/4}\cap\mathcal F^{\le 1/4}$ contains many  intriguing $3$-graphs. For instance, it is easy to check that $K_4^{(3)-} \in \mathcal F^{\ge 1/4}\cap\mathcal F^{\le 1/4}$ (for more details, refer to the proof of Theorem~\ref{7-k43-}). 
	Recall that $F^{\star}_5$ is a $3$-graph
	that is obtained from $K^{(3)-}_4$
	by adding a new vertex whose link graph forms a matching
	on the vertices of $K^{(3)-}_4$ (see Figure~\ref{fig:picture1}).
	From that point of view, the shadow $\partial K_4^{(3)-}$ of $K_4^{(3)-}$ can be described as the complete graph $K_4$ that can be partitioned into exactly three
	matchings of size $2$. Given $t\in \{5, 6, 7\}$, let $F^{\star}_t$ denote the $3$-graph on $t$ vertices that is obtained from $K_4^{(3)-}$ by adding  $(t-4)$ vertices such that the link graph of each additional vertex forms a matching of size 2 on $V(K_4^{(3)-})$ and the link graphs of additional vertices are disjoint (see Figure~\ref{fig:3-pictures}). Furthermore, we consider a $3$-graph, denoted by $\hat F^{\star}_7$, obtained by adding a hyperedge consisting of three additional vertices to the $3$-graph $F^{\star}_7$.
	In other words, $\hat F^{\star}_7$ is formed by $7$ vertices $\{a,b,c,d,x,y,z\}$ and $10$ hyperedges $\{abc, abd, acd, xab, xcd, yac, ybd, zad, zbc, xyz\}$. Using Theorem~\ref{main-thm1}, we can determine the uniform Tur\'an density of $\hat F^{\star}_7$.
	
	\begin{theorem}\label{7-k43-}
		We have $\pi_{\points}(\hat F^{\star}_7)=1/4$.
	\end{theorem}

	\begin{remark}
		Observe that $K^{(3)-}_4\sub F^{\star}_5\sub F^{\star}_6 \sub F^{\star}_7 \sub \hat F^{\star}_7$, but $\hat F^{\star}_7$ is not a blow-up of $K^{(3)-}_4$, $F^{\star}_5$, $F^{\star}_6$ or $F^{\star}_7$.
		In particular, due to $K^{(3)-}_4\sub F^{\star}_5\sub F^{\star}_6 \sub F^{\star}_7 \sub \hat F^{\star}_7$, Theorem~\ref{7-k43-} implies that $\{F^{\star}_5, F^{\star}_6, F^{\star}_7,\hat F^{\star}_7\} \subset\mathcal F^{\ge 1/4}\cap\mathcal F^{\le 1/4}$.
	\end{remark}

	We next present more examples of $3$-graphs $F$ with $F\in \mathcal F^{\ge 1/4}\cap\mathcal F^{\le 1/4}$.
	Let $t \geq 4$ be an integer and $G$ be a graph on $(t-1)$ vertices with chromatic number $\chi(G) = 3$. Suppose $G^a_t$ is the 3-graph constructed from a vertex $a$ and the graph $G$ as follows: $V(G^a_t) = \{a\} \cup V(G)$ and $E(G^a_t) = \{\{a,x,y\} : \{x,y\} \in E(G)\}$. Furthermore, let $\mathcal G^{(3)}_t$ denote the family consisting of all those 3-graphs $G^a_t$.
	Using Theorem~\ref{main-thm1}, we can also readily obtain the following result.
	
	\begin{theorem}\label{thm:3-chromatic}
		For any integer $t \geq 4$ and any 3-graph $G^a_t\in \mathcal G^{(3)}_t$, we have $\pi_{\points}(G^a_t)=1/4$.
	\end{theorem}

	\begin{remark}
		Observe that each 3-graph $G^a_t\in \mathcal G^{(3)}_t$ is a subhypergraph of the blow-up of $K^{(3)-}_4$. Therefore, we can also determine $\pi_{\points}(G^a_t)$ by supersaturation.
		Furthermore, Theorem~\ref{thm:3-chromatic} can be proved using the method described by Reiher, R\"{o}dl, and Schacht~\cite{k43minus-2}, although it was not explicitly stated in their work.
	\end{remark}

	Finally we notice an interesting class of 3-graphs.
	Given $t\ge 3$, 
	let $W^{(3)}_t$ denote the \emph{$3$-uniform wheel}  on $t$ vertices, which is a $3$-graph with $t$ vertices that can be ordered $\{v_0, v_1, \dots, v_{t-1}\}$ such that  $\{v_0, v_i, v_{i+1}\}$ is a $3$-edge of $W^{(3)}_t$ for $i\in \mathbb Z_{t}\setminus\{0\}$ (see Figure~\ref{figure-wheel}). Clearly, the shadow $\partial W^{(3)}_t$ of $W^{(3)}_t$ forms a wheel graph on $t$ vertices. In addition, it is easy to  check that $W^{(3)}_t\in G^{(3)}_t$ when $t$ is even; $W^{(3)}_t$  satisfies the property of Theorem~\ref{0-density} when $t$ is odd (because it is a $3$-partite $3$-graphs). Therefore, by Theorem~\ref{thm:3-chromatic} and Theorem~\ref{0-density}, we can directly determine the uniform Tur\'an density of $W^{(3)}_t$ for all $t\ge 3$.

	\begin{figure}
		\centering
		\begin{subfigure}[b]{0.45\textwidth}
			\centering
			\begin{tikzpicture}
				[inner sep=2pt,
				vertex/.style={circle, draw=blue!50, fill=blue!50},
				]
				\filldraw[black] (6,3.5) circle (2pt);
				\node at (5.7,3.6) {$v_0$};
				\filldraw[black] (3,2.2) circle (2pt);
				\node at (2.8,2) {$v_1$};
				\filldraw[black] (4,2.2) circle (2pt);
				\node at (3.8,2) {$v_2$}; 
				\filldraw[black] (5,2.2) circle (2pt);
				\filldraw[black] (5.5,2.2) circle (2pt);
				\draw[line width=1pt,black,opacity=1] (5,2.2) -- (5.5,2.2);
				\filldraw[black] (6,2.2) circle (1pt);
				\filldraw[black] (6.1,2.2) circle (1pt);
				\filldraw[black] (6.2,2.2) circle (1pt);
				\filldraw[black] (6.5,2.2) circle (2pt);
				\filldraw[black] (7,2.2) circle (2pt);
				\draw[line width=1pt,black,opacity=1] (6.5,2.2)-- (7,2.2);
				\filldraw[black] (7.8,2.2) circle (2pt);
				\node at (8,1.9) {$v_{t-1}$}; 
				\filldraw[black] (9,2.2) circle (2pt);
				\node at (9.2,1.9) {$v_{1}$}; 
				\draw[line width=1pt,black,opacity=1] (6,3.5) -- (3,2.2);
				\draw[line width=1pt,black,opacity=1] (6,3.5) -- (4,2.2);
				\draw[line width=1pt,black,opacity=1] (3,2.2)-- (4,2.2);
				\draw[line width=1pt,black,opacity=1] (7.8,2.2)--  (9,2.2);
				\draw[line width=1pt,black,opacity=1] (6,3.5) -- (7.8,2.2);
				\draw[line width=1pt,black,opacity=1] (6,3.5)--  (9,2.2);
				\draw[line width=1pt,black,opacity=1] (6,3.5) -- (5,2.2);
				\draw[line width=1pt,black,opacity=1] (6,3.5)--  (5.5,2.2);
				\draw[line width=1pt,black,opacity=1] (6,3.5) -- (6.5,2.2);
				\draw[line width=1pt,black,opacity=1] (6,3.5)--  (7,2.2);
				\filldraw[black] (7.2,2.2) circle (1pt);
				\filldraw[black] (7.3,2.2) circle (1pt);
				\filldraw[black] (7.4,2.2) circle (1pt);
				\filldraw[black] (4.4,2.2) circle (1pt);
				\filldraw[black] (4.5,2.2) circle (1pt);
				\filldraw[black] (4.6,2.2) circle (1pt);
			\end{tikzpicture}
			\caption{The $3$-uniform wheel  $W^{(3)}_t$.}\label{figure-wheel}
		\end{subfigure}
		\hfill
		\begin{subfigure}[b]{0.5\textwidth}
			\centering
			\begin{tikzpicture}
				[inner sep=2pt,
				vertex/.style={circle, draw=blue!50, fill=blue!50},
				]
				\filldraw[black] (6,2.5) circle (2pt);
				\node at (5.8,2.6) {$x$};
				\filldraw[black] (3,1.2) circle (2pt);
				\node at (2.8,1) {$v_1$};
				\filldraw[black] (4,1.2) circle (2pt);
				\node at (3.8,1) {$v_2$}; 
				\filldraw[black] (5,1.2) circle (2pt);
				\filldraw[black] (5.5,1.2) circle (2pt);
				\draw[line width=1pt,black,opacity=1] (5,1.2) -- (5.5,1.2);
				\filldraw[black] (6,1.2) circle (1pt);
				\filldraw[black] (6.1,1.2) circle (1pt);
				\filldraw[black] (6.2,1.2) circle (1pt);
				\filldraw[black] (6.5,1.2) circle (2pt);
				\filldraw[black] (7,1.2) circle (2pt);
				\draw[line width=1pt,black,opacity=1]  (6.5,1.2) --  (7,1.2);
				\filldraw[black] (7.8,1.2) circle (2pt);
				\node at (8,0.9) {$v_{t-1}$}; 
				\filldraw[black] (9,1.2) circle (2pt);
				\node at (9.2,0.9) {$v_{1}$}; 
				\filldraw[black] (6,-0.3) circle (2pt);
				\node at (5.8,-0.5) {$y$}; 
				\draw[line width=1pt,black,opacity=1] (6,2.5) -- (3,1.2);
				\draw[line width=1pt,black,opacity=1] (6,2.5) -- (4,1.2);
				\draw[line width=1pt,black,opacity=1] (3,1.2)-- (4,1.2);
				\draw[line width=1pt,black,opacity=1] (6,-0.3) -- (3,1.2);
				\draw[line width=1pt,black,opacity=1] (6,-0.3)-- (4,1.2);
				\draw[line width=1pt,black,opacity=1] (6,-0.3) -- (7.8,1.2);
				\draw[line width=1pt,black,opacity=1] (6,-0.3)--  (9,1.2);
				\draw[line width=1pt,black,opacity=1] (7.8,1.2)--  (9,1.2);
				\draw[line width=1pt,black,opacity=1] (6,2.5) -- (7.8,1.2);
				\draw[line width=1pt,black,opacity=1] (6,2.5)--  (9,1.2);
				\draw[line width=1pt,black,opacity=1] (6,-0.3) -- (5,1.2);
				\draw[line width=1pt,black,opacity=1] (6,-0.3)--  (5.5,1.2);
				\draw[line width=1pt,black,opacity=1] (6,2.5) -- (5,1.2);
				\draw[line width=1pt,black,opacity=1] (6,2.5)--  (5.5,1.2);
				\draw[line width=1pt,black,opacity=1] (6,-0.3) -- (6.5,1.2);
				\draw[line width=1pt,black,opacity=1] (6,-0.3)--  (7,1.2);
				\draw[line width=1pt,black,opacity=1] (6,2.5) -- (6.5,1.2);
				\draw[line width=1pt,black,opacity=1] (6,2.5)--  (7,1.2);
				\filldraw[black] (7.2,1.2) circle (1pt);
				\filldraw[black] (7.3,1.2) circle (1pt);
				\filldraw[black] (7.4,1.2) circle (1pt);
				\filldraw[black] (4.4,1.2) circle (1pt);
				\filldraw[black] (4.5,1.2) circle (1pt);
				\filldraw[black] (4.6,1.2) circle (1pt);
			\end{tikzpicture}
			\caption{A double pyramid on $t+1$ vertices.}\label{figure-double-pyramid}
		\end{subfigure}
		\caption{}
	\end{figure}

	\begin{corollary}\label{wheel}
		For  $t\ge 3$,
		\[
		\pi_{\points}(W^{(3)}_t) = \begin{cases}
			1/4, &\text{~if~} t \text{ is even};\\
			0, &\text{~if~} t \text{ is odd}.
		\end{cases}
		\]
	\end{corollary}

	\begin{remark}
		A  \emph{double pyramid} is a $3$-graph on $t+ 1$ vertices $x, y, v_1,\dots, v_{t-1}$ for
		some integer $t\ge 4$, whose hyperedges are $\{x, v_i, v_{i+1}\}$ and $\{y, v_i, v_{i+1}\}$ for $i\in \mathbb Z_{t}\setminus\{0\}$,  see Figure~\eqref{figure-double-pyramid}. The double pyramids have very good connectivity and symmetry. In the study of $3$-graphs, the double pyramid has been of interest because of its connections to other structures and properties, such as Brown, Erd\H{o}s and S\'os~\cite{double-pyramid-1973} proved that if a $3$-graph $H$ has $n$ vertices and $\Omega(n^{5/2})$ edges, then $H$ contains a double pyramid. 
		By Corollary~\ref{wheel}, we can directly obtain the uniform Tur\'an density of all double pyramids with a given size, since every double pyramid on $t+ 1$ vertices is a blow-up of $W^{(3)}_t$.  Note that we have that $\pi(F)\ge \pi_{\points}(F)$ for every $3$-graph $F$. Therefore, the  Tur\'an density of a double pyramid with odd vertices is at least $1/4$.
	\end{remark}

	\subsection{Organization}
	The rest of this paper is organized as follows. In the next section, we will review a framework presented in Reiher's survey~\cite{reiher2020extremal}, namely, \emph{reduced hypergraphs} (or \emph{reduced $3$-graphs}). This framework encompasses regularity arguments for $3$-graphs and has been a widely recognized tool for embedding subhypergraphs in dense $3$-graphs.
	In Section~\ref{sec-Inter}, we will prove some intersection lemmas for reduced $3$-graphs that match the scenarios used in our arguments. We will then apply these lemmas to analyze  the structural properties of reduced $3$-graphs with a density greater than $1/4$ in Section~\ref{sec-struc}, and give an embedding lemma (see Lemma~\ref{embed-lem}), which is the key to proving Theorem~\ref{main-thm1}. At the end of Section~\ref{sec-struc}, we will give the proof of Theorem~\ref{main-thm1}.
	In Section~\ref{sec-pf-less-than1/4}, we will  prove  some auxiliary lemmas;
	then use them to give a proof of Theorem~\ref{less-than1/4}.
	Finally, using Theorem~\ref{main-thm1}, we will prove Theorem~\ref{7-k43-} and Theorem~\ref{thm:3-chromatic} in Section~\ref{Applications}.

	\section{Reduced hypergraphs}\label{sec-Reduced}
	In this section, we introduce the concept of reduced $3$-graphs, which encapsulates hypergraph
	regularity arguments about the uniform Tur\'an density of $3$-graphs. We follow some notions introduced in Reiher's survey~\cite[Section 3]{reiher2020extremal}.
	It should be noted that the concept of  reduced $3$-graphs has also been mentioned in other related works~\cite{k43minus-2,vanishing,cycle-uniform,1/27,chen2022beyond}.

	\begin{definition}[Reduced $3$-graphs] \label{def-reduced-graph}
		Given a finite index set $I\subset \mathbb N$, for each pair $ij\in [I]^{2}$, let $\mathcal {P}_{ij}$ denote a finite nonempty vertex set  such that for any two distinct  pairs $ij, i'j'\in [I]^{2}$ the sets $\mathcal {P}_{ij}$ and $\mathcal {P}_{i'j'}$ are disjoint.
		For any triple $ijk \in [I]^{3}$, let $\mathcal A_{ijk}$ denote a $3$-partite 
		$3$-graph with vertex partition $\mathcal {P}_{ij}$, $\mathcal {P}_{jk}$ and $\mathcal {P}_{ik}$.
		Then the $\binom{|I|}{2}$-partite 
		$3$-graph $\mathcal A$ with
		\[
		V(\mathcal A)=\bigcup_{ij\in [I]^{2}}\mathcal {P}_{ij} ~~~~\text{and}~~~~ E(\mathcal A)=\bigcup_{ijk\in [I]^{3}} E(\mathcal A_{ijk})
		\]
		is called an {\it $I$-reduced} $3$-graph with \emph{index set} $I$, \emph{vertex classes} $\mathcal {P}_{ij}$  and {\it constituents} $\mathcal A_{ijk}$. 
		Moreover,  we refer  to vertices that belong to $\mathcal {P}_{ij}$ as {\it left vertices}
		of $\mathcal A_{ijk}$, those that belong to $\mathcal {P}_{jk}$ as {\it right vertices} of $\mathcal A_{ijk}$ and those that belong to $\mathcal {P}_{ik}$ as {\it top vertices} of $\mathcal A_{ijk}$.
	\end{definition}
	
\begin{remark}\label{rmk-reverse}
For brevity, we often simply write ``$\mathcal A$ is a reduced $3$-graph" instead of ``$\mathcal A$ is an $I$-reduced $3$-graph for some index set $I \subset \mathbb N$''.
Given an $I$-reduced $3$-graph $\mathcal A$ with $I=x_1x_2\dots x_n$,
		the \emph{reverse} $\mathcal A^{-1}$  of $\mathcal A$ is a reduced $3$-graph with the same vertex set and the same hyperedge set as $\mathcal A$
		but with the partition of vertices given by $\mathcal P^{-1}_{x_i,x_j}=\mathcal P_{x_{n-j+1},x_{n-i+1}}$ for $ij\in [n]^2$.
		Clearly,  every reduced $3$-graph $\mathcal A$ and its reverse $\mathcal A^{-1}$ are isomorphic. In particular, for a triple $ijk\in [n]^3$, left vertices with respect to $\mathcal A_{ijk}$ become right vertices of $\mathcal A^{-1}_{ijk}$; right vertices with respect to $\mathcal A_{ijk}$ becomes left vertices of $\mathcal A^{-1}_{ijk}$; top vertices with respect to $\mathcal A_{ijk}$ are still top for $\mathcal A^{-1}_{ijk}$. 
	\end{remark}
	
	Similar to the role of reduced graphs in Szemer\'edi's regularity method, we need to express our density conditions in terms of reduced $3$-graphs.

	\begin{definition}
		Given $d \in[0, 1] $ and an $I$-reduced $3$-graph $\mathcal A$, we say that  $\mathcal A$ is {\it $d$-dense} if 
		\[
		|E(\mathcal A_{ijk})|\ge d |\mathcal P_{ij}||\mathcal P_{jk}||\mathcal P_{ik}|
		\]
		holds for all $ijk\in [I]^{3}$.
	\end{definition}

	When applying the regularity method for 3-graphs, if a given 3-graph $F$ can be embedded into a uniformly dense $3$-graph $H$, then the corresponding reduced $3$-graph $\mathcal A$ with respect to $H$ inherits a property reflecting this.
	As a result, we need to define how a reduced $3$-graph $\mathcal{A}$ ``embeds" a given $3$-graph $F$.

	\begin{definition}\label{def-embed}
		Given a $3$-graph $F$ and an  $I$-reduced $3$-graph  $\mathcal A$, we say that $\mathcal A$ {\it embeds} $F$ if there is a pair of maps $(\phi,\psi)$ such that
		\begin{itemize} 
			\item[$\rm (1)$]  $\phi:V(F)\to I$ and $\psi: \partial F\to V(\mathcal A)$;
			\item[$\rm (2)$] if $\{u,v\}\in \partial F$, then $\psi(u, v) \in \mathcal P_{\phi (u)\phi(v)}$;
			\item[$\rm (3)$] if $\{u,v,w\} \in E(F)$, then $\{\psi(u,v), \psi(v,w), \psi(u,w)\}\in E(\mathcal A_{\phi(u)\phi(v)\phi(w)})$.
		\end{itemize}
	\end{definition}
	
	In Reiher's survey~\cite{reiher2020extremal}, a general theorem is provided that reduces  the process of proving the upper bound of the uniform Tur\'an density of a $3$-graph $F$ to embedding $F$ into reduced $3$-graphs with the same density.

	\begin{theorem}[Reiher~\cite{reiher2020extremal}]\label{thm-turan-reduced}
		Let $F$ be a $3$-graph and $d\in [0,1]$.
		If for any $\eps>0$ there exists $N\in \mathbb N$ such that each $(d+\eps)$-dense $[N]$-reduced $3$-graph $\mathcal A$ embeds $F$, then $\pi_{\points}(F)\le d$.
	\end{theorem}

	\section{Intersection lemmas}\label{sec-Inter}
	In this section,  we present some intersection lemmas for reduced $3$-graphs that match the scenarios used in our arguments. 
	To prove these lemmas, a main tool is the classical Ramsey theorem for multicolored
	$k$-graphs, which is stated below for reference.

	\begin{theorem}[Ramsey~\cite{Ramsey}]  \label{thm-Ramsey} 
		For any $r_R, k_R, n_R \in \mathbb N$, there exists $N_R \in \mathbb N$ such that every $r_R$-edge-coloring of a $k_R$-uniform clique with $N_R$ vertices contains a monochromatic $k_R$-uniform clique with $n_R$ vertices.
	\end{theorem}

	Given an $[N]$-reduced $3$-graph  $\mathcal A$, with candidate sets of good properties for each constituent, such as considering the candidate sets as top vertices with large degrees for each constituent, we aim to choose a representative vertex from each vertex class in $\mathcal A$ that possesses good properties for each constituent that it belongs to. Since each vertex class in $\mathcal A$ may belong to many different constituents, it is possible that no single vertex is suitable for all constituents of $\mathcal A$, even if the size of the candidate sets with good properties relative to each constituent involving the vertex class is linearly proportional. However, by leveraging the power of Ramsey theory, it is possible to find such a representative vertex by passing to the induced subhypergraph of $\mathcal A$.

	The next two lemmas are intended to identify representative vertices based on candidate sets from the top vertices and the right  vertices of the constituents, respectively.
	It should be noted that the proofs of these two lemmas have been implicitly presented in~\cite{vanishing} using an iterative approach, while alternative proofs can also be found in~\cite{1/27, cycle-uniform}.

	\begin{lemma}[{\cite[Lemma 4.3]{cycle-uniform}} or {\cite[Lemma 10]{1/27}}]\label{top-lem}
		For any $\delta>0$ and $n\in \mathbb N$, there exists  $N\in \mathbb N$ such that the following holds.
		If $\mathcal A$ is an $[N]$-reduced $3$-graph and there exists a subset $X_{ijk}\subseteq \mathcal P_{ik}$ with $|X_{ijk}|\ge \delta |\mathcal P_{ik}|$ for all $ijk \in [N]^3$, then there is an induced subhypergraph $\mathcal A'\subseteq \mathcal A$ with index set $I$ of size $n$ and there are vertices $x_{ik}\in \mathcal P_{ik}$ for all $ik\in [I]^2$ such that
		\[
		x_{ik}\in \bigcap_{i<j<k,~j\in I}  X_{ijk}.
		\]
	\end{lemma}

	\begin{lemma}[{\cite[Lemma 4.5]{cycle-uniform}} or {\cite[Lemma 9]{1/27}}]\label{right-lem-1}
		For any $\delta>0$  and $n\in \mathbb N$, there exists  $N\in \mathbb N$ such that the following holds.
		If $\mathcal A$ is an $[N]$-reduced $3$-graph and there exists a subset $X_{ijk}\subseteq \mathcal P_{jk}$ with $|X_{ijk}|\ge \delta |\mathcal P_{jk}|$ for all $ijk \in [N]^3$, then there is an induced subhypergraph $\mathcal A' \subseteq \mathcal A$ with index set $I$ of size $n$ and vertices $x_{jk}\in \mathcal P_{jk}$ for all $jk\in [I]^2$ such that 
		\[
		x_{jk}\in \bigcap_{i<j<k,~i\in I}  X_{ijk}.
		\]
	\end{lemma}

	The following intersection lemma is designed to select representative vertices from candidate sets based on  $3$-partite graphs that we define later on sets $\mathcal P_{ij}\cup \mathcal P_{ij'} \cup \mathcal P_{ik}$ for indices $i<j<j'<k$. Given an $[N]$-reduced $3$-graph  $\mathcal A$, for any $ijj'k\in [N]^4$, if we have a candidate set $R_{ijj'k}\subseteq \mathcal P_{ij'}$ that is linearly proportional in size, then by Ramsey theory, we can still find a representative vertex with respect to those sets $P_{ij'}$ by passing to the induced subhypergraph of $\mathcal A$. In fact, fixed $i$ and $j'$, the candidate indicators of $j$ and $k$ can form a dense bipartite graph, and any bipartite graph with positive density contains a complete bipartite graph of a given size.
	Note that the following lemma is equivalent to~\cite[Lemma 4.5]{cycle-uniform}, and as a result, we omit its proof here.

	\begin{lemma}[{\cite[Lemma 4.8]{cycle-uniform}}]\label{right-lem}
		For any $\delta>0$ and $n\in \mathbb N$, there exists  $N\in \mathbb N$ such that the following holds.
		If $\mathcal A$ is an $[N]$-reduced $3$-graph and there exists a subset $R_{ijj'k}\subseteq \mathcal P_{ij'}$ with $|R_{ijj'k}|\ge \delta |\mathcal P_{ij'}|$ for each $ijj'k\in [N]^4$, then there is an induced subhypergraph $\mathcal A'\subseteq \mathcal A$ with index set $I$ of size $n$ and there are vertices $r_{ij'}\in \mathcal P_{ij'}$ for all  $ij'\in [I]^2$ such that
		\[
		r_{ij'}\in \bigcap_{i<j<j'<k,~jk\in [I]^2}  R_{ijj'k}.
		\]
	\end{lemma}

	Given an $[N]$-reduced $3$-graph  $\mathcal A$,  if we have a dense
	candidate set $L_{ijj'k}\subseteq \mathcal P_{ij}$ for any $ijj'k\in [N]^4$, we can not obtain the analogous conclusion as in Lemma~\ref{right-lem}. It is possible that there may not be any vertex $\ell_{ij}\in \mathcal P_{ij}$ such that $\ell_{ij}\in L_{ijxy}\cap L_{ijxz}\cap L_{ijyz}$ for any $1\le i<j<x<y<z\le N$. In fact,  for fixed $i$ and $j$, the candidate indicators of $j'$ and $k$ form a dense graph, but a graph with positive density does not necessarily contain a complete graph of a given size; it may not even contain a triangle. 
	However, we may still need to identify the representative vertex from some vertex class of $\mathcal A$ that match the scenarios used in our later arguments. Therefore, we give the following conclusion.

	\begin{lemma}\label{left-lem}
		For any $\delta>0$ and $n_1,n_2\in \mathbb N$, there exists  $N\in \mathbb N$ such that the following holds.
		If $\mathcal A$ is an $[N]$-reduced $3$-graph and there exists a subset $L_{ijj'k}\subseteq \mathcal P_{ij}$ with $|L_{ijj'k}|\ge \delta |\mathcal P_{ij}|$ for  each $ijj'k\in [N]^4$, then there is an induced subhypergraph $\mathcal A'\subseteq \mathcal A$ with index set $I_1\cup I_2$ where $|I_1|=n_1$, $|I_2|=n_2$ and $\max I_1<\min I_2$, and there exist vertices  $\ell_{ij}\in \mathcal P_{ij}$ for all $ij\in [I_1]^2$  such that 
		\[
		\ell_{ij}\in \bigcap_{i<j<j'<k,~ j'\in I_1,~ k\in I_2}  L_{ijj'k}.
		\]
	\end{lemma}

	\begin{proof}
		Given $\delta>0$ and $n_1,n_2\in \mathbb N$, we apply Theorem~\ref{thm-Ramsey} with  $r_R= 2$, $k_R=n_1+n_2+2$ and $n_R =\max\{2n_1+n_2,  n_0+3\}  $ to get $N\in \mathbb N$, where $n_0$ is a positive integer such that any balanced bipartite graph
		on $n_0$ vertices with at least $\frac{\delta n_0^2}{4}$ edges contains a complete bipartite graph $K_{n_1,n_2}$. 
		Let $\mathcal A$ be an $[N]$-reduced $3$-graph satisfying the condition of the lemma.
		We now consider a $2$-edge-colored $(n_1+n_2+2)$-uniform clique on the vertex set $[N]$ as follows.
		Given an $(n_1+n_2+2)$-edge $Q= a_1a_2\dots a_{n_1+n_2+2} \in [N]^{n_1+n_2+2}$, we define that $Q$ is colored blue if there exists a vertex  $\ell_{a_1a_2}\in \mathcal P_{a_1a_2}$ such that 
		\[
		\ell_{a_1a_2} \in \bigcap_{2<i\le n_1+2< j\le n_1+n_2+2,~a_ia_j\in [Q]^2} L_{a_1a_2a_ia_j},
		\]
		otherwise $Q$ is colored red.	
		
		By Theorem~\ref{thm-Ramsey}, there exists a subset $S\subseteq [N]$ with $|S|=n_R$ such that all edges $Q$ induced on $S$ have same color. For convenience, we rearrange the indices in $S$ and write $S=[n_R]$. 
		Note that each $L_{12xy}\subseteq \mathcal P_{12}$ satisfies
		$|L_{12xy}|\ge \delta |\mathcal P_{12}|$ for all  $2<x\le \lceil \frac{n_R}{2} \rceil+1 <y\le n_R$. 
		By double counting, there exists a vertex $\ell_{12}\in \mathcal P_{12}$ that  is contained in at least   
		\[
		\delta\left \lceil \frac{n_R-2}{2} \right\rceil \left\lfloor \frac{n_R-2}{2} \right\rfloor >
		\frac{\delta n_0^2}{4}
		\] 
		$L_{12xy}$'s for  some $x\in \{3,\dots, \lceil \frac{n_R}{2} \rceil+1\}$ and $y\in \{\lceil \frac{n_R}{2} \rceil+2,\dots, n_R\}$. 
		Let us define a bipartite graph $G$ with partition sets $\{3,\dots, \lceil \frac{n_R}{2} \rceil+1\}$ and $\{\lceil \frac{n_R}{2} \rceil+2,\dots, n_R\}$ by putting an edge $\{x,y\}\in E(G)$ if and only if $\ell_{12}\in   L_{12xy}$. Clearly, $|E(G)|>
		\frac{\delta n_0^2}{4}$.
		Based on the choice of $n_0$, there exists a complete bipartite graph $K_{n_1,n_2}$ in $G$, which means that there is an $(n_1+n_2+2)$-set $Q= 12a_3\dots a_{n_1+n_2+2}$ with $\{a_3, \dots, a_{n_1+2 }\}\subset \{3,\dots, \lceil \frac{n_R}{2} \rceil+1 \}$ and $\{a_{n_1+3}, \dots,   a_{n_1+n_2+2}\}\subset \{\lceil \frac{n_R}{2} \rceil+2, \dots,  n_R\}$ such that $Q$ is colored blue. Therefore, the common color for the edges induced on $S$ is  blue.
		
		Now we choose $I_1=\{1, 2,  \dots, n_1\}$, $I_2=\{n_R-n_2+1, \dots, n_R\}$. Clearly, $I_1\cup I_2 \subset [n_R]$, $|I_1|=n_1$, $|I_2|=n_2$  and $\max I_1<\min I_2$.
		For any $ij\in [I_1]^2$,
		we extend  $\{i,j\}\cup \{x\in I_1:x>j\} \cup I_2$ 
		to an $(n_1+n_2+2)$-set $Q$ by adding some elements from set $\{n_1+1,\dots, 2n_1 \}$. Since $Q$ is colored blue,  there is a vertex  $\ell_{ij}\in \mathcal P_{ij}$ satisfying
		\[
		\ell_{ij}\in \bigcap_{i<j<j'<k,~ j'\in I_1,~ k\in I_2}  L_{ijj'k}.
		\]
	\end{proof}

	In the proofs of the next section, we will encounter a situation that appears unfavorable. 
	Given an $[N]$-reduced $3$-graph $\mathcal A$, we may only have sets of candidate properties for some constituents of $\mathcal A$, but not all. 
	For example, only consider candidate sets for those $\mathcal A_{ijk}$ with $i<j<t$ and $t<k$ for some  $t\in [N]$. Obviously, we can not directly use Ramsey's theorem to reach the analogous conclusion as in Lemma~\ref{top-lem} or Lemma~\ref{right-lem-1}, since Ramsey's theorem is based on complete hypergraphs. 
	The following intersection lemma is designed to select representative vertices  based on candidate sets depending on  some, but not all, constituents of $\mathcal A$.
	In this case, it is feasible to select the representative vertices by using Ramsey's theorem iteratively.

	\begin{lemma}\label{bi-right-lem}
		For any $\delta>0$ and $n_1, n_2\in \mathbb N$, there exists  $N\in \mathbb N$ such that the following holds.
		If $\mathcal A$ is an $[N]$-reduced $3$-graph and there exist subsets $W_{ijk}\subseteq \mathcal P_{jk}$ with $|W_{ijk}|\ge \delta |\mathcal P_{jk}|$ for each $1\le i<j\le N-n_2<k\le N$, then there is an induced subhypergraph  $\mathcal A'\subseteq \mathcal A$ with index set $I\cup \{N-n_2+1,\dots, N\}$ where $I \subseteq [N-n_2]$ and $|I|=n_1$, and there exist vertices  $w_{jk}\in \mathcal P_{jk}$ for all $j\in I$ and $k\in \{N-n_2+1,\dots, N\}$  such that 
		\[
		w_{jk}\in \bigcap_{i<j<k,~i\in I} W_{ijk}.
		\]
		
	\end{lemma}
	
	\begin{proof}
		Given $\delta>0$ and  $n_1, n_2\in \mathbb N$, we apply Theorem~\ref{thm-Ramsey} iteratively $n_2$ times, and choose constants satisfying the following hierarchy (from right to left):
		\[
		N-n_2=m_{n_2}\gg t_{n_2}\gg m_{n_2-1}\gg \dots \gg t_2
		\gg m_1\gg t_1\gg m_0=n_1,1/\delta.
		\]
		To be more precise,  in the first iteration, we apply Theorem~\ref{thm-Ramsey} with $r_R=2$, $k_R=m_0$ and $n_R=t_1=\max\{2n_1, \lceil \frac{n_1}{\delta}\rceil \}$ to get $N_R=m_{1}$.
		In the $\ell$th  iteration, we apply Theorem~\ref{thm-Ramsey} with $r_R=2$, $k_R=m_{\ell-1}$ and $n_R=t_\ell=\max\{2m_{\ell-1}, \lceil \frac{m_{\ell-1}}{\delta}\rceil \}$ to get $N_R=m_{\ell}$.  In the last iteration, we apply Theorem~\ref{thm-Ramsey} with $r_R=2$, $k_R=m_{n_2-1}$ and $n_R=t_{n_2}=\max\{2m_{n_2-1}, \lceil \frac{m_{n_2-1}}{\delta}\rceil \}$ to get $N_R=m_{n_2}$, and set $N=m_{n_2}+n_2$.

		Let $\mathcal A$ be an $[N]$-reduced $3$-graph satisfying the condition of the lemma. First, let $k=N-n_2+1$. 
		We consider a $2$-edge-coloring $m_{n_2-1}$-uniform clique on the vertex set $[m_{n_2}]$ as follows. Given an $m_{n_2-1}$-edge $Q=a_1a_2\dots a_{m_{n_2-1}} \in [m_{n_2}]^{m_{n_2-1}}$, let $Q$ be colored blue if there exists a vertex  $w_{a_{m_{n_2-1}}k}\in \mathcal P_{a_{m_{n_2-1}}k}$ such that 
		\[
		w_{a_{m_{n_2-1}}k} \in \bigcap_{i\in Q\setminus \{a_{m_{n_2-1}}\}} W_{ia_{m_{n_2-1}}k},
		\]
		otherwise $Q$ is colored red.	
		
		By Theorem~\ref{thm-Ramsey}, there exists a set $T_{n_2}\subseteq [m_{n_2}]$ with $|T_{n_2}|=t_{n_2}$ such that all hyperedges $Q$ induced on set $T_{n_2}$ have same color. For convenience, we rearrange the indices in $T_{n_2}$ and write $T_{n_2}=[t_{n_2}]$.
		By the condition of the lemma, each  $W_{it_{n_2}k}\subseteq \mathcal P_{t_{n_2}k}$ satisfies
		$|W_{it_{n_2}k}|\ge \delta |\mathcal P_{t_{n_2}k}|$ for all  $i\in [t_{n_2}-1]$.
		By double counting, there exists a vertex $w_{t_{n_2}k}\in \mathcal P_{t_{n_2}k}$ that  is contained in at least  
		\[
		\frac{\delta|\mathcal P_{t_{n_2}k}|(t_{n_2}-1)}{|\mathcal P_{t_{n_2}k}|}=\delta (t_{n_2}-1)  \ge m_{n_2-1}-1
		\] 
		$W_{it_{n_2}k}$'s for some $i\in [t_{n_2}-1]$, where the last inequality follows from $t_{n_2}=\max\{2m_{n_2-1}, \lceil \frac{m_{n_2-1}}{\delta}\rceil \}$. This implies  that
		the common color for the hyperedges induced on $[t_{n_2}]$ is blue.
		Now we choose $M_{n_2-1}=\{m_{n_2-1}+1, \dots, 2m_{n_2-1}\}$. For any $j\in M_{n_2-1}$, consider the $m_{n_2-1}$-set $Q=\{j-m_{n_2-1}+1,\dots,j\}$. Since $Q$ is colored blue,  there is a vertex  $w_{jk}\in \mathcal P_{jk}$ satisfying
		\[
		w_{jk}\in \bigcap_{i\in Q\setminus \{j\}} W_{ijk} \subseteq \bigcap_{i<j<k,~i\in M_{n_2-1}}  W_{ijk}.
		\]

		Next, let $K_\ell=\{N-n_2+1,\dots, N-n_2+\ell\}$ for $\ell\in [n_2]$. For convenience, we rearrange the indices in $M_{n_2-\ell+1}$ and write $M_{n_2-\ell+1}=[m_{n_2-\ell+1}]$.
		We apply Theorem~\ref{thm-Ramsey} again, but this time to consider a $2$-edge-coloring $m_{n_2-\ell}$-uniform clique on the vertex set $[m_{n_2-\ell+1}]$ as follows. Given an $m_{n_2-\ell}$-edge $Q=a_1a_2\dots a_{m_{n_2-\ell}} \in [m_{n_2-\ell+1}]^{m_{n_2-\ell}}$, we define that $Q$ is colored blue if there exists a vertex  $w_{a_{m_{n_2-\ell}}k}\in \mathcal P_{a_{m_{n_2-\ell}}k}$ for $k\in K_\ell\setminus K_{\ell-1}$  (let $K_0=\emptyset$) such that 
		\[
		w_{a_{m_{n_2-\ell}}k} \in \bigcap_{i\in Q\setminus \{a_{m_{n_2-\ell}}\}} W_{ia_{m_{n_2-\ell}}k},
		\]
		otherwise $Q$ is colored red.	
		Similar to the first step, there exists a set $T_{n_2-\ell+1}\subseteq [m_{n_2-\ell+1}]$ with $|T_{n_2-\ell+1}|=t_{n_2-\ell+1}$ such that all hyperedges $Q$ induced on set $T_{n_2-\ell+1}$ is colored blue. Therefore, we can choose a set $M_{n_2-\ell}\subset T_{n_2-\ell+1}$ such that for each $j\in M_{n_2-\ell}$, there is a vertex $w_{jk}\in \mathcal P_{jk}$ with $k\in K_\ell\setminus K_{\ell-1}$ satisfying
		\[
		w_{jk}\in \bigcap_{i<j<k,~i\in M_{n_2-\ell}}  W_{ijk}.
		\]
		
		After performing the iteration procedure $n_2$ times as described, we obtain sets  $K_{n_2}=\{N-n_2+1,\dots, N\}$, $I\subset T_1\subset M_1\subset \dots \subset  T_{n_2}\subset [m_{n_2}]$ with $|I|=n_1$, and  vertices  $w_{jk}\in \mathcal P_{jk}$ for all $j\in I$ and $k\in K_{n_2}$  such that 
		\[
		w_{jk}\in \bigcap_{i<j<k,~i\in I}  W_{ijk}.
		\]
	\end{proof}

	\section{Structural results}\label{sec-struc}
	In this section, we will extract some good properties from the reduced $3$-graphs $\mathcal A$ with a density greater than $1/4$ by passing to the induced subhypergraph of $\mathcal A$, which guarantees the existence of  certain structures of  $\mathcal A$.

	We first project each constituent of reduced $3$-graphs onto an appropriate bipartite graph, similar to the approach used in~\cite{k43minus-2,chen2022beyond}.
	Given a $(1/4+\eps)$-dense $I$-reduced $3$-graph $\mathcal A$ (without loss of generality, let $\eps\ll 1$),  for any $ijk\in [I]^3$
	we define a bipartite graph $Q^i_{jk}$ with bipartition $\{\mathcal P_{ij}, \mathcal P_{ik}\}$ by putting an edge between
	$u\in \mathcal P_{ij}$ and $v\in \mathcal P_{ik}$ if and only if there exist at least $\eps^2|\mathcal P_{jk}|$ vertices $w\in \mathcal P_{jk}$  such that $\{u,v,w\}\in E(\mathcal A_{ijk})$. 
	Similarly, let $Q^k_{ij}$ denote the bipartite graph on the vertex set $\{\mathcal P_{ik}, \mathcal P_{jk}\}$ by putting an edge between
	$u\in \mathcal P_{ik}$ and $v\in \mathcal P_{jk}$ if and only if there exist at least $\eps^2|\mathcal P_{ij}|$ vertices $w\in \mathcal P_{ij}$  such that $\{u,v,w\}\in E(\mathcal A_{ijk})$. We write $d_{Q^i_{jk}}(v)$ and $d_{Q^k_{ij}}(v)$  for the degree of $v \in  V(Q^i_{jk})$ or $v \in  V(Q^k_{ij})$, respectively. Based on the above definitions and Cauchy-Schwartz inequality, Chen and Sch{\"u}lke~\cite{chen2022beyond} gave the following observation.

	\begin{observation}[\cite{chen2022beyond}]  \label{degree-square} 
		For any $\eps>0$ and $n\in \mathbb N$, there exists  $N\in \mathbb N$ such that the following holds.
		For each $(1/4+\eps)$-dense $[N]$-reduced $3$-graph $\mathcal A$, there exists $\mathcal A'\subseteq \mathcal A$ induced on set $I\subseteq [N]$ of size $n$ such that at least one of the following holds:
		\stepcounter{propcounter}
		\begin{enumerate}[label = ({\bfseries \Alph{propcounter}\arabic{enumi}})]
			\item \label{p-a} For all  $ijk\in [I]^3$,  $\sum_{v\in \mathcal P_{ik}} d_{Q^i_{jk}}(v)^2 \ge (\frac{1}{4}+\frac{\eps}{2})|\mathcal P_{ij}|^2|\mathcal P_{ik}|$. 
			\item \label{p-b} For all  $ijk\in [I]^3$,  $\sum_{v\in \mathcal P_{ik}} d_{Q^k_{ij}}(v)^2 \ge (\frac{1}{4}+\frac{\eps}{2})|\mathcal P_{jk}|^2|\mathcal P_{ik}|$.
		\end{enumerate}
	\end{observation}

	Given an $[n]$-reduced $3$-graph $\mathcal A$, if $\mathcal A$ satisfies the property~\ref{p-b} in Observation~\ref{degree-square}, clearly by Remark~\ref{rmk-reverse}, the reverse $\mathcal A^{-1}$ of $\mathcal A$ satisfies the property~\ref{p-a} in Observation~\ref{degree-square}, and vice versa. In addition, since every reduced $3$-graph $\mathcal A$ and its reverse $\mathcal A^{-1}$ are completely isomorphic, we directly obtain the following observation.

	\begin{observation}\label{degree-square-2} 
		For any $\eps>0$ and $n\in \mathbb N$, there exists  $N\in \mathbb N$ such that the following holds.
		For each $(1/4+\eps)$-dense $[N]$-reduced $3$-graph $\mathcal A$, there exists an  induced subhypergraph $\mathcal A'\subseteq \mathcal A$ with index set $I$ of size $n$ such that for all  $ijk\in [I]^3$,  
		\begin{equation}\label{eq-a}
			\sum_{v\in \mathcal P_{ik}} d_{Q^i_{jk}}(v)^2 \ge (\frac{1}{4}+\frac{\eps}{2})|\mathcal P_{ij}|^2|\mathcal P_{ik}|.
		\end{equation}
	\end{observation}
	
	Next, we give a definition to collect these top vertices with large degree in $Q^i_{jk}$. Given a $(1/4+\eps)$-dense $I$-reduced $3$-graph $\mathcal A$,  for each $ijk\in [I]^3$ and positive integer $r$, let
	\[
	S^j_{ik}(r):=\{v\in \mathcal P_{ik}: d_{Q^i_{jk}}(v)\ge (1/2+r\eps^2 )|\mathcal P_{ij}|\}.
	\]
	
	Using Observation~\ref{degree-square-2} and Theorem~\ref{thm-Ramsey}, we have the following result. We note that the proof of Lemma~\ref{S-r} is implied in the statement in~\cite[Section 2]{k43minus-2} and~\cite[Section 3]{chen2022beyond}, and thus, we only include its sketch for completeness.

	\begin{lemma}\label{S-r} 
		For any $\eps>0$ and $n\in \mathbb N$, there exist positive integers $N$ and $r_*$ such that the following holds.
		Suppose that $\mathcal A$ is an $[N]$-reduced $3$-graph satisfying Inequality~\eqref{eq-a} for all $ijk\in [N]^3$. Then there exists an induced subhypergraph $\mathcal A'\subseteq \mathcal A$ on set $I$ of size $n$ such that for all  $ijk\in [I]^3$ we have
		\begin{equation}\label{ineq1}
			\sum_{v\in \mathcal P_{ik}} d_{Q^i_{jk}}(v)^2 \ge (\frac{1}{4}+\frac{\eps}{2})|\mathcal P_{ij}|^2|\mathcal P_{ik}|,~~~
			|S^j_{ik}(r_*)| \ge \eps^2|\mathcal P_{ik} | ~~~\text{and}~~~|S^j_{ik}(r_*+1)| < \eps^2|\mathcal P_{ik}|.
		\end{equation}
	\end{lemma}
	
	\begin{proof}
		We apply Theorem~\ref{thm-Ramsey} with  $r_R= \lfloor (2\eps^2)^{-1} \rfloor$, $k_R=3$
		and $n_R =n $ to get $N\in \mathbb N$. Let $\mathcal A$ be an $[N]$-reduced $3$-graph  satisfying Inequality~\eqref{eq-a} for all $ijk\in [N]^3$. Due to $\mathcal A$ satisfying Inequality~\eqref{eq-a} for all $ijk\in [N]^3$,  we have
		\[
		(\frac{1}{4}+\frac{\eps}{2})|\mathcal P_{ij}|^2|\mathcal P_{ik}|\le \sum_{v\in \mathcal P_{ik}} d_{Q^i_{jk}}(v)^2 \le    |S^j_{ik}(1)||\mathcal P_{ij}|^2+ |\mathcal P_{ik}\setminus S^j_{ik}(1)|\left((1/2+\eps^2)|\mathcal P_{ij}|\right)^2,
		\]
		and we have $|S^j_{ik}(1)|\ge \eps^2|\mathcal P_{ik} |$  by $\eps\ll 1$. 
		Moreover, $S^i_{jk}(r)=\emptyset$ for $r>(2\eps^2)^{-1}$. Thus, we may color each triple $ijk\in [N]^3$ by
		the maximum $r\in [r_R]$ for which $|S^j_{ik}(r)|\ge \eps^2|\mathcal P_{ik} |$. By Theorem~\ref{thm-Ramsey}, there exists a subset $I\subseteq [N]$ of size $n$ and a positive integer $r_*$ such that  the desired conclusion holds.
	\end{proof}

	Now we restrict ourselves to the $I$-reduced $3$-graphs satisfying Inequalities~\eqref{ineq1}. For any $i\in I$, let $Q^i:=\bigcup_{ijk\in [I]^3} Q^i_{jk}$. Chen and Sch{\"u}lke~\cite{chen2022beyond} proved the following result, which is a robust version of a result in~\cite{k43minus-2}.

	\begin{lemma} [{\cite[Lemma 3.2]{chen2022beyond}}]\label{many-triangles} 
		Let $\mathcal A$ be an $I$-reduced $3$-graph satisfying Inequalities~\rm{(\ref{ineq1})}. Then for any $ijj'k\in [I]^4$ and $v\in S^{j}_{ik}(r_*)\cap S^{j'}_{ik}(r_*)$, there are at least $\eps^2 |\mathcal P_{ij}||\mathcal P_{ij'}|$ triangles $\{v,u, w\}$ in $Q^i$ with $u\in \mathcal P_{ij}$ and $w \in \mathcal P_{ij'}$.
	\end{lemma}
	
	We are now ready to prove an embedding lemma of reduced $3$-graphs with density greater than $1/4$, which is the main result of
	this section. The lemma will be used to upper bound the uniform Tur\'an density
	of $3$-graphs $F$ with $F\in \mathcal F^{\le 1/4}$.
	
	\begin{lemma}\label{embed-lem} 
		For any $\eps>0$ and $m_1, m_2\in \mathbb N$, there exists  $N\in \mathbb N$ such that each $(1/4+\eps)$-dense $[N]$-reduced $3$-graph $\mathcal A$ contains an induced subhypergraph $\mathcal A'\sub \mathcal A$ with index set $M_1\cup M_2$ such that $|M_1|=m_1$, $|M_2|=m_2$, $\max M_1<\min M_2$ and $\mathcal A'$ satisfies the following properties:
		\begin{itemize} 
			\item[$\rm (1)$]   
			There exist vertices $x_{ij}, y_{ij}, w^1_{ij}\in \mathcal P_{ij}$ for all $ij\in [M_1\cup M_2]^2$ such that $\{y_{ij}, w^1_{jk}, x_{ik}\}\in E(\mathcal A_{ijk})$ for all $ijk\in [M_1\cup M_2]^3$.
			\item[$\rm (2)$] There also exist vertices $z_{ij}, w^2_{ij}, w^3_{ij} \in \mathcal P_{ij}$ for all $ij\in [M_1]^2$ such that 
			$\{z_{ij}, w^2_{jk}, y_{ik}\}\in E(\mathcal A_{ijk})$ for all $ijk\in [M_1]^3$, and 	$\{z_{ij}, w^3_{jk}, x_{ik}\}\in E(\mathcal A_{ijk})$ for all $ij\in [M_1]^2$ and $k\in M_2$.
		\end{itemize}
	\end{lemma}

	\begin{proof}
		Given $\eps>0$ and  $m_1, m_2\in \mathbb N$, we choose constants satisfying the following hierarchy (from right to left):
		\[
		N\gg N_7\gg N_6\gg N_5\gg N_4 \gg N_3
		\gg N_2\gg N_1\gg m_1, m_2, 1/\eps.
		\]
		To be more precise, we determine $N$ as follows: first apply Lemma~\ref{bi-right-lem} with $\delta=\eps^2$, $n_1=m_1$ and $n_2=m_2$ to get $N_1$, then Lemma~\ref{right-lem-1} with $\delta=\eps^2$ and $n=N_1-m_2$ to get $N_2$, then Lemma~\ref{left-lem} with $\delta=\eps^2/2$, $n_1=N_2$ and $n_2=m_2$ to get $N_3$, then apply Lemma~\ref{right-lem-1} again with $\delta=\eps^2$ and $n=N_3$ to get $N_4$, then Lemma~\ref{right-lem} with $\delta=\eps^2/2$ and $n=N_4$ to get $N_5$, then 
		Lemma~\ref{top-lem} with $\delta=\eps^2$  and $n=N_5$ to get $N_6$, then Lemma~\ref{S-r} with $\eps$ and $n=N_6$ to get $N_7$ and $r_*$, and finally apply Observation~\ref{degree-square-2} with $\eps$ and $n=N_7$ to get $N$.
		Let  $\mathcal A$ be a $(1/4+\eps)$-dense $[N]$-reduced $3$-graph. By Observation~\ref{degree-square-2}, there exists an induced subhypergraph $\mathcal A_7\subseteq \mathcal A$ with index set $I_7$ of size $N_7$ that satisfies Inequality~\eqref{eq-a} for all $ijk\in [I_7]^3$.
		
		For the sake of clarity, the rest of the proof will be processed by the following six steps.
		
		\noindent{\bf Step 1: choose the vertices $x_{ij}$.} We apply Lemma~\ref{S-r} to $\mathcal A_7$, and obtain an induced subhypergraph $\mathcal A_6\subseteq \mathcal A_7$ with index set $I_6$ of size $N_6$ satisfying Inequalities~(\ref{ineq1}). In particular, for all $ijk\in [I_6]^3$ we have
		\[ 
		S^j_{ik}(r_*)\subseteq \mathcal P_{ik} ~~~\text{and}~~~|S^j_{ik}(r_*)| \ge \eps^2|\mathcal P_{ik} | .
		\]
		We now apply Lemma~\ref{top-lem} to $\mathcal A_6$ and sets $S^j_{ik}(r_*)$, obtaining an induced subhypergraph $\mathcal A_5\subseteq \mathcal A_6$ on  set $I_5$ of size $N_5$, and vertices $x_{ik}\in \mathcal P_{ik}$ for all $ik\in [I_5]^2$ such that 
		\[ 
		x_{ik}\in \bigcap_{i<j<k,~j\in I_5} S^j_{ik}(r_*).
		\]
		
		\noindent{\bf Step 2: choose the vertices $y_{ij}$.} For each $ijj'k\in [I_5]^4$, we define $R_{ijj'k}\subseteq \mathcal P_{ij'}$ to be the set  of vertices $v$ in $\mathcal P_{ij'}$ such that $\{v, x_{ik}\}$ is contained in at least $\eps^2/2 |\mathcal P_{ij}|$ triangles in the graph $Q^i[\mathcal P_{ij}\cup \mathcal P_{ij'}\cup \mathcal P_{ik}]$. Since $x_{ik}\in  S^j_{ik}(r_*)\cap S^{j'}_{ik}(r_*)$, by Lemma~\ref{many-triangles}, there are at least $\eps^2|\mathcal P_{ij}||\mathcal P_{ij'}|$ edges in the graph $Q^i[\mathcal P_{ij}\cup \mathcal P_{ij'}]$. By the definition of $R_{ijj'k}$, we obtain that
		\[
		\eps^2|\mathcal P_{ij}||\mathcal P_{ij'}|\le |R_{ijj'k}||\mathcal P_{ij}|+ \eps^2/2 |\mathcal P_{ij}| |\mathcal P_{ij'} \setminus R_{ijj'k}|,
		\]
		which implies that $|R_{ijj'k}|\ge \frac{\eps^2}{2} |\mathcal P_{ij'}|$.
		By Lemma~\ref{right-lem} applied with $\mathcal A_5$ and sets $R_{ijj'k}$, there exists an induced subhypergraph $\mathcal A_4\subseteq \mathcal A_5$ with index set $I_4$ of size $N_4$, and vertices $y_{ij'}\in \mathcal P_{ij'}$ for all $ij'\in [I_4]^2$ such that 
		\[ 
		y_{ij'}\in \bigcap_{i<j<j'<k,~jk\in [I_4]^2} R_{ijj'k}.
		\]
		
		\noindent{\bf Step 3: choose the vertices $w^1_{ij}$.} For each $ijk\in [I_4]^3$, let 
		\[
		W^1_{ijk}=\left\{w\in \mathcal P_{jk}: \{y_{ij}, w, x_{ik}\}\in E(\mathcal A_{ijk})\right\}.
		\]
		Since $\{y_{ij}, x_{ik}\}\in E(Q^i_{jk})$, by the definition of $Q^i_{jk}$, we have  $|W^1_{ijk}|\ge \eps^2|\mathcal P_{jk}|$ for all $ijk\in [I_4]^3$. Apply  Lemma~\ref{right-lem-1} to $\mathcal A_4$ and sets $W^1_{ijk}$, to get an induced subhypergraph $\mathcal A_3\subseteq \mathcal A_4$ with index set $I_3$ of size $N_3$, and vertices $w^1_{jk}\in \mathcal P_{jk}$ for all $jk\in [I_3]^2$ such that 
		\[ 
		w^1_{jk}\in \bigcap_{i<j<k,~i\in I_3} W^1_{ijk}.
		\]
		Until now, we have $x_{ij}, y_{ij}, w^1_{ij}\in \mathcal P_{ij}$ for all $ij\in [I_3]^2$ and $\{y_{ij}, w^1_{jk}, x_{ik}\}\in E(\mathcal A_{ijk})$ for all $ijk\in [I_3]^3$.
		
		\noindent{\bf Step 4: choose the vertices $z_{ij}$.} For any $ijj'k\in [I_3]^4$,
		Let 
		\[
		L_{ijj'k}=\{z\in \mathcal P_{ij}: \{z, y_{ij'}, x_{ik}\} \text{~forms~a~triangle~in~} Q^i[\mathcal P_{ij}\cup \mathcal P_{ij'}\cup \mathcal P_{ik}] \}.
		\]
		Due to the choice of $y_{ij'}$ in Step 2, we have
		$|L_{ijj'k}|\ge \eps^2/2 |\mathcal P_{ij}|$. By Lemma~\ref{left-lem} applied with $\mathcal A_3$ and sets $L_{ijj'k}$, there exists an induced subhypergraph $\mathcal A_2\subseteq \mathcal A_3$ on set $I_2\cup M_2$ with $|I_2|=N_2$, $|M_2|=m_2$ and $\max I_2< \min M_2$,  and there exist vertices $z_{ij}\in \mathcal P_{ij}$ for all $ij\in [I_2]^2$ such that 
		\[ 
		z_{ij}\in \bigcap_{i<j<j'<k,~j'\in I_2,~ k\in M_2} L_{ijj'k}.
		\]
		
		\noindent{\bf Step 5: choose the vertices $w^2_{ij}$.} 
		For each $ijk\in [I_2]^3$, let us define
		\[
		W^2_{ijk}=\left\{w\in \mathcal P_{jk}: \{z_{ij}, w, y_{ik}\}\in E(\mathcal A_{ijk})\right\}.
		\]
		Similar to Step 3, we have  $|W^2_{ijk}|\ge \eps^2|\mathcal P_{ik}|$ since $\{z_{ij}, y_{ik}\}\in E(Q^i_{jk})$.  We apply  Lemma~\ref{right-lem-1} with the induced subhypergraph $\mathcal A'_2\subset \mathcal A_2$ on set $I_2$ and the sets $W^2_{ijk}$, to get an induced subhypergraph $\mathcal A'_1\subseteq \mathcal A'_2$ on set $I_1\subseteq I_2$ of size $N_1$, and vertices $w^2_{jk}\in \mathcal P_{jk}$ for all $jk\in [I_1]^2$ such that
		\[ 
		w^2_{jk}\in \bigcap_{i<j<k,~i\in I_1} W^2_{ijk}.
		\]
		
		\noindent{\bf Step 6: choose the vertices $w^3_{ij}$.} Finally let $\mathcal A_1$ be the reduced $3$-graph induced with index set $I_1\cup M_2$, and define
		\[
		W^3_{ijk}=\left\{w\in \mathcal P_{jk}: \{y_{ij}, w, x_{ik}\}\in E(\mathcal A_{ijk})\right\}
		\]
		for all $ij\in [I_1]^2$ and $k\in M_2$.
		Since $\{y_{ij}, x_{ik}\}\in E(Q^i_{jk})$,
		we have$|W^3_{ijk}|\ge \eps^2|\mathcal P_{jk}|$.
		Then we apply Lemma~\ref{bi-right-lem} with $\mathcal A_1$ and the sets $W^3_{ijk}$  to obtain an induced subhypergraph $\mathcal A'\subseteq \mathcal A_1$ on set $M_1\cup M_2$ with $M_1 \subseteq [I_1]$ and $|M_1|=m_1$, and vertices $w^3_{jk}\in \mathcal P_{jk}$ for all $j\in M_1$ and $k\in M_2$ such that
		\[
		w^3_{jk}\in \bigcap_{i<j<k,~i\in M_1} W^3_{ijk}.
		\]
		
		Hence, the reduced $3$-graph  $\mathcal A'$ with index set $M_1\cup M_2$ and the vertices $x_{ij}, y_{ij}, z_{ij}$, $w^1_{ij}$, $w^2_{ij}$ and $w^3_{ij}$  satisfy the conclusion of the lemma.
	\end{proof}

	Now, we will prove our main theorem using Lemma~\ref{embed-lem}.
	\begin{proof}[{\bf Proof of Theorem~\ref{main-thm1}}]
		Fix a $3$-graph $F\in \mathcal F^{\ge 1/4}\cap \mathcal F^{\le 1/4}$ with $f$ vertices. By Fact~\ref{low-bound}, we have $\pi_{\points}(F)\ge 1/4$.
		Next, by Theorem~\ref{thm-turan-reduced}, we only need to show
		that for every $\eps>0$, there exists $N\in \mathbb N$ such that every $(1/4+\eps)$-dense $[N]$-reduced $3$-graph $\mathcal A$ can  embed $F$.
		
		Due to $F\in \mathcal F^{\le 1/4}$, there exists an ordering $(v_1, \dots, v_f)$ of $V(F)$ with an integer $i^\star\in [f]$, and
		there is a  $6$-coloring 
		$
		\chi: \partial F\to \{\textcolor{red}{red},\textcolor{blue}{blue}, \textcolor{green}{green},\textcolor{violet}{violet},\textcolor{cyan}{cyan}, \textcolor{black}{black}\}
		$
		such that every $hyperedge$ $\{v_i, v_j, v_k\}$ of $F$ with $i<j< k$ satisfies exactly one of the following:
		\begin{enumerate}
			\item $\left(\chi(v_i,v_j), \chi(v_j,v_k), \chi(v_i,v_k)\right)=(\textcolor{blue}{blue}, \textcolor{violet}{violet},\textcolor{red}{red})$ and $ijk\in [f]^3$;
			\item  $\left(\chi(v_i,v_j), \chi(v_j,v_k), \chi(v_i,v_k)\right)=(\textcolor{green}{green}, \textcolor{cyan}{cyan},\textcolor{blue}{blue})$ and $ijk\in [i^\star]^3$;
			\item  $\left(\chi(v_i,v_j), \chi(v_j,v_k), \chi(v_i,v_k)\right)=(\textcolor{green}{green}, \textcolor{black}{black},\textcolor{red}{red})$ and $1\le i<j<i^\star< k\le f$.
		\end{enumerate}
		We apply Lemma~\ref{embed-lem} with $\eps>0$, $m_1=i^\star$ and $m_2=f-i^\star$ to get
		$N$. Let $\mathcal A$ be a $(1/4+\eps)$-dense  $[N]$-reduced $3$-graph. By Lemma~\ref{embed-lem}, there exists an induced subhypergraph $\mathcal A'\subseteq \mathcal A$ 
		with index set $M_1\cup M_2$ satisfying the properties given in Lemma~\ref{embed-lem}, where $|M_1|=m_1$, $|M_2|=m_2$ and $\max M_1<\min M_2$.
		
		For brevity, let $M_1=\{1, \dots, {m_1}\}$, $M_2=\{{m_1+1}, \dots, f\}$  and $\phi:V(F)\to M_1\cup M_2$ satisfying $\phi(v_\ell)=\ell$ for all $\ell \in [f]$. For any pair $\{v_i,v_j\}\in \partial F$ with $ij\in [f]^2$, we choose the vertex $x_{ij}$ from $\mathcal P_{ij}$ if $\chi (v_i,v_j)$ is \textcolor{red}{red}, choose the vertex $y_{ij}$ from $\mathcal P_{ij}$ if $\chi (v_i,v_j)$ is \textcolor{blue}{blue}, choose the vertex $z_{ij}$ from $\mathcal P_{ij}$ if $\chi (v_i,v_j)$ is \textcolor{green}{green}, choose the vertex $w^1_{ij}$ from $\mathcal P_{ij}$ if $\chi (v_i,v_j)$ is \textcolor{violet}{violet}, choose the vertex $w^2_{ij}$ from $\mathcal P_{ij}$ if $\chi (v_i,v_j)$ is \textcolor{cyan}{cyan},
		and choose the vertex $w^3_{ij}$ from $\mathcal P_{ij}$ if $\chi (v_i,v_j)$ is \textcolor{black}{black}. Recalling Definition~\ref{def-embed}, the induced subhypergraph $\mathcal A'$ can embed exactly $F$.
	\end{proof}
	
	\section{The proof of Theorem~\ref{less-than1/4}}\label{sec-pf-less-than1/4}
	In this section, we will prove Theorem~\ref{less-than1/4}, and the proof of Theorem~\ref{less-than1/4} is based on Theorem~\ref{0-density} and the following lemma.
	
	\begin{lemma}\label{less-than1/4-lem}
		Let $F$ be a $3$-graph  such that the following holds.
		There exists an  order $(v_1, \dots, v_{|V(F)|})$ of $V(F)$ and a  $5$-coloring $\chi: \partial F\to \{ \textcolor{red}{red}, {\color{blue}blue}, {\color{green}green},\textcolor{violet}{violet},\textcolor{cyan}{cyan} \}$  
		such that every $3$-edge $\{v_i, v_j, v_k\}$ of $F$ with $i<j< k$ satisfies
		\[
		\left(\chi(v_i,v_j), \chi(v_j,v_k), \chi(v_i,v_k)\right)=(\textcolor{blue}{blue}, \textcolor{violet}{violet},\textcolor{red}{red}) \text{~or~}  (\textcolor{green}{green}, \textcolor{cyan}{cyan},\textcolor{blue}{blue}).
		\]
		Then $\pi_{\points}(F)\le d^*$.
	\end{lemma}

	To prove Lemma~\ref{less-than1/4-lem}, we need some auxiliary lemmas. Therefore, we first prove Theorem~\ref{less-than1/4}, and postpone the proof of Lemma~\ref{less-than1/4-lem} to the end of this section.
	
	\begin{proof}[{\bf Proof of Theorem~\ref{less-than1/4}}] Given a $3$-graph $F\in \mathcal F^{\le 1/4}$, let $\sigma=(v_1, \dots, v_{|V(F)|})$  be an order of $V(F)$ and $\chi: \partial F\to \{\textcolor{red}{red},\textcolor{blue}{blue}, \textcolor{green}{green},\textcolor{violet}{violet},\textcolor{cyan}{cyan}, \textcolor{black}{black}\}$ be a $6$-coloring provided by property {\rm ($\spadesuit$)}. 
		
		If $E_1=\emptyset$ in~\ref{p11}, then we will construct a $3$-coloring 
		$\chi_1: \partial F\to \{ \textcolor{green}{green},\textcolor{cyan}{cyan}, \textcolor{blue}{blue}\}$ which is obtained from the $6$-coloring $\chi$ by replacing \textcolor{violet}{violet}, \textcolor{cyan}{cyan} and \textcolor{black}{black} with \textcolor{cyan}{cyan}, and \textcolor{red}{red} and \textcolor{blue}{blue}  with \textcolor{blue}{blue}.
		Observe that  every hyperedge $\{v_i, v_j, v_k\}$ of $F$ with $i<j< k$ satisfies
		\[ 
		\left(\chi(v_i,v_j), \chi(v_j,v_k), \chi(v_iv_k)\right)=(\textcolor{green}{green}, \textcolor{cyan}{cyan},\textcolor{blue}{blue}),
		\]
		under the order $\sigma$ and $3$-coloring $\chi_1$. By Theorem~\ref{0-density}, we have $\pi_{\points}(F)=0<d^*$.
		
		If $E_2=\emptyset$ in~\ref{p22}, then we consider a $3$-coloring 
		$\chi_2: \partial F\to \{\textcolor{blue}{blue}, \textcolor{violet}{violet}, \textcolor{red}{red} \}$, which is also obtained from the $6$-coloring $\chi$ by replacing \textcolor{violet}{violet}, \textcolor{cyan}{cyan} and \textcolor{black}{black} with \textcolor{violet}{violet}, and \textcolor{green}{green} and \textcolor{blue}{blue}  with \textcolor{blue}{blue}.
		Observe that  every hyperedge $\{v_i, v_j, v_k\}$ of $F$ with $i<j< k$ satisfies
		\[ 
		\left(\chi(v_i,v_j), \chi(v_j,v_k), \chi(v_iv_k)\right)=(\textcolor{blue}{blue}, \textcolor{violet}{violet}, \textcolor{red}{red}),
		\]
		under the order $\sigma$ and $3$-coloring $\chi_2$. By Theorem~\ref{0-density}, $\pi_{\points}(F)=0<d^*$.
		
		Finally, we consider $E_3=\emptyset$ in~\ref{p33}. Since $E_3=\emptyset$, there is no pair $\{v_i,v_j\}\in \partial F$ such that $\chi(v_i,v_j)=\textcolor{black}{black}$.
		Clearly,  $F$ with the ordering $\sigma$ and coloring $\chi$ satisfies the properties given in Lemma~\ref{less-than1/4-lem}. Hence, $\pi_{\points}(F)\le d^*$.
	\end{proof}
	
	Next, we state and prove several lemmas that are useful for the proof of Lemma~\ref{less-than1/4-lem}. 
	The first lemma follows from Lemma~\ref{right-lem-1} by applying it to the reverse of reduced $3$-graph $\mathcal A$ given in the statement of Lemma~\ref{right-lem-1}. We note that it also appears in~\cite{cycle-uniform,1/27}.
	
	\begin{lemma}[{\cite[Lemma 4.4]{cycle-uniform}} or {\cite[Lemma 8]{1/27}}]\label{left-lem-1}
		For any $\delta>0$  and $n\in \mathbb N$, there exists  $N\in \mathbb N$ such that the following holds.
		If $\mathcal A$ is an $[N]$-reduced $3$-graph and there exists a subset $X_{ijk}\subseteq \mathcal P_{ij}$ with $|X_{ijk}|\ge \delta |\mathcal P_{ij}|$ for each $ijk\in [N]^3$, then there is an $\mathcal A' \subseteq \mathcal A$ induced on set $I$ of size $n$ and there are vertices $x_{ij}\in \mathcal P_{ij}$ for all $ij\in [I]^2$ such that 
		\[
		x_{ij}\in \bigcap_{i<j<k,~k\in I}  X_{ijk}.
		\]
	\end{lemma}

	Given an $I$-reduced $3$-graph $\mathcal A$, $ijk\in [I]^3$ and $\delta>0$, we define
	\begin{itemize} 
		\item $L^k_{ij}:=\{v \in \mathcal P_{ij}:  v \text{~is~contained~in~at~least~} \delta|\mathcal P_{jk}||\mathcal P_{ik}| \text{~hyperedges~in~}  \mathcal A_{ijk}\}$, 
		\item $R^i_{jk}:=\{v \in \mathcal P_{jk}:  v \text{~is~contained~in~at~least~} \delta|\mathcal P_{ij}||\mathcal P_{ik}| \text{~hyperedges~in~}  \mathcal A_{ijk}\}$, and
		\item $T^j_{ik}:=\{v \in \mathcal P_{ik}:  v \text{~is~contained~in~at~least~} \delta|\mathcal P_{ij}||\mathcal P_{jk}| \text{~hyperedges~in~}  \mathcal A_{ijk}\}$,
	\end{itemize} 
	
	\begin{lemma}[{\cite[Lemma 6]{1/27}}]\label{almost-same}
		For any $\delta>0$ and $n\in \mathbb N$, there exists $N \in \mathbb N$ such that the following holds. For every $[N]$-reduced $3$-graph $\mathcal A$, there exist reals $t_i$ for $i\in \{1,2,3\}$ and an induced subhypergraph $\mathcal A'\subseteq \mathcal A$ with index set $I$ of size $n$ such that for all $ijk\in [I]^3$, we have
		\begin{equation}\label{eq-almost-same}
			t_1\le |L^k_{ij}|/|\mathcal P_{ij}| \le t_1+\delta,~~ t_2\le |R^i_{jk}|/|\mathcal P_{jk}| \le t_2+\delta  ~~\text{and}~~ t_3\le |T^j_{ik}| /|\mathcal P_{ik}|\le t_3+\delta.
		\end{equation}
	\end{lemma}
	
	\begin{lemma}\label{left-top-right}
		For every $\eps>0$, there exists $\delta>0$ such that for every $n\in \mathbb N$ there exists $N\in \mathbb N$ such that the following holds. For every $(d^*+\eps)$-dense $[N]$-reduced $3$-graph $\mathcal A$, there exists an induced subhypergraph $\mathcal A' \subseteq \mathcal A$ with index set $I$ of size $n$ such that $|L^{\ell}_{ik}\cap T^{j}_{ik}|\ge \delta|\mathcal P_{ik}|$ for all  $ijk\ell \in [I]^4$.
	\end{lemma}
	
	\begin{proof} Given $\eps>0$, we choose $\delta=\eps^3/10$.
		We first apply Theorem~\ref{thm-Ramsey} with  $r_R= 2$, $k_R=5$ and $n_R =2n+1$ to get $N_1\in \mathbb N$, and then apply Lemma~\ref{almost-same}  with $\delta$ and $N_1$ to get $N$. Let $\mathcal A$ be a $(d^*+\eps)$-dense $[N]$-reduced $3$-graph.
		Using Lemma~\ref{almost-same} with $\delta$ and $N_1$, there exists an induced subhypergraph $\mathcal A_1\subseteq \mathcal A$ on set $I_1$ of size $N_1$ and 
		reals $t_i$ for $i\in \{1,2,3\}$ satisfying Inequalities~\eqref{eq-almost-same}. For convenience, let $L:=L^k_{ij}$, $R:=R^i_{jk}$ and $T:=T^j_{ik}$. We first claim that 
		\begin{equation}\label{eq-t_3}
			t_3\ge d^*+\eps/2.
		\end{equation}
		Indeed,  for every $ijk\in [I_1]^3$, we have 
		\[
		(d^*+\eps)|\mathcal P_{ij}||\mathcal P_{jk}||\mathcal P_{ik}|\le|E(\mathcal A_{ijk})| \le |T||\mathcal P_{ij}||\mathcal P_{jk}|+ \delta|\mathcal P_{ij}||\mathcal P_{jk}||\mathcal P_{ik}\setminus T|
		\le (t_3+2\delta)|\mathcal P_{ij}||\mathcal P_{jk}||\mathcal P_{ik}|,
		\]
		which means that $t_3\ge d^*+\eps/2$ since $\delta<\eps/10$. 
		Further, we have 
		\begin{equation}\label{eq-t123}
			t_1+t_2+t_3\ge  3 \sqrt[3]{d^*}+\eps/3.
		\end{equation}
		If not, suppose that $t_1+t_2+t_3< 3 \sqrt[3]{d^*}+\eps/3$. By  Lemma~\ref{almost-same}, for every $ijk\in [I_1]^3$ we have 
		\[
		\begin{split}
			|E(\mathcal A_{ijk})| & \le \delta|\mathcal P_{jk}||\mathcal P_{ik}||\mathcal P_{ij}\setminus L|+
			\delta|\mathcal P_{ij}||\mathcal P_{ik}||\mathcal P_{jk}\setminus R|+
			\delta|\mathcal P_{ij}||\mathcal P_{jk}| |\mathcal P_{ik}\setminus T|+ |L||R||T|\\
			&\le 3\delta |\mathcal P_{ij}||\mathcal P_{jk}||\mathcal P_{ik}|+ (t_1+\delta)(t_2+\delta)(t_3+\delta)|\mathcal P_{ij}||\mathcal P_{jk}||\mathcal P_{ik}|\\
			&\le \left(3\delta+ \left(\frac{t_1+t_2+t_3+3\delta}{3}\right)^3\right) |\mathcal P_{ij}||\mathcal P_{jk}||\mathcal P_{ik}|< (d^*+\eps/2) |\mathcal P_{ij}||\mathcal P_{jk}||\mathcal P_{ik}|,
		\end{split}
		\]
		which contradicts that $\mathcal A$ is $(d^*+\eps)$-dense, where the last inequality follows from $\delta<\eps/10$.

		Next, we consider a $2$-edge-colored $5$-uniform clique on set $I_1$ as follows. Given a $5$-edge $Q=a_1a_2a_3a_4a_5 \in [I_1]^{5}$, if $|L^{a_5}_{a_2a_4}\cap T^{a_3}_{a_2a_4}|\ge \frac{\eps}{4}|\mathcal P_{a_2a_4}|$, $Q$ is colored red; otherwise, if $|R^{a_1}_{a_2a_4}\cap T^{a_3}_{a_2a_4}|\ge \frac{\eps}{4}|\mathcal P_{a_2a_4}|$, $Q$ is colored blue. If neither of these two cases occurs, that is, $|L^{a_5}_{a_2a_4}\cap T^{a_3}_{a_2a_4}|< \frac{\eps}{4}|\mathcal P_{a_2a_4}|$ and $|R^{a_1}_{a_2a_4}\cap T^{a_3}_{a_2a_4}|< \frac{\eps}{4}|\mathcal P_{a_2a_4}|$  which means that
		\[
		|(L^{a_5}_{a_2a_4} \cup R^{a_1}_{a_2a_4}) \cap T^{a_3}_{a_2a_4}|=|(L^{a_5}_{a_2a_4}\cap T^{a_3}_{a_2a_4}) \cup ( R^{a_1}_{a_2a_4}\cap T^{a_3}_{a_2a_4})|< \frac{\eps}{2}|\mathcal P_{a_2a_4}|.
		\]
		Further, we have
		\[\begin{split}
			|L^{a_5}_{a_2a_4}\cup R^{a_1}_{a_2a_4}\cup T^{a_3}_{a_2a_4}| &=|L^{a_5}_{a_2a_4}\cup R^{a_1}_{a_2a_4}|+|T^{a_3}_{a_2a_4}|-|(L^{a_5}_{a_2a_4} \cup R^{a_1}_{a_2a_4}) \cap T^{a_3}_{a_2a_4}|\\
			&= |L^{a_5}_{a_2a_4}|+|R^{a_1}_{a_2a_4}|+|T^{a_3}_{a_2a_4}|-|L^{a_5}_{a_2a_4}\cap R^{a_1}_{a_2a_4}|-|(L^{a_5}_{a_2a_4} \cup R^{a_1}_{a_2a_4}) \cap T^{a_3}_{a_2a_4}|\\
			&> (t_1+t_2+t_3)|\mathcal P_{a_2a_4}|-\left( |\mathcal P_{a_2a_4}\setminus T^{a_3}_{a_2a_4}|+\frac{\eps}{4}|\mathcal P_{a_2a_4}| \right)-\frac{\eps}{2}|\mathcal P_{a_2a_4}|\\
			&\ge \left( (t_1+t_2+t_3)-1+t_3-\eps/4-\eps/2\right)|\mathcal P_{a_2a_4} |\\
			&\stackrel{\eqref{eq-t_3}, \eqref{eq-t123}}{> }\left(3 \sqrt[3]{d^*}+d^*-1+\eps/13\right)|\mathcal P_{a_2a_4}|,
		\end{split}
		\]
		which is impossible since $L^{a_5}_{a_2a_4}\cup R^{a_1}_{a_2a_4}\cup T^{a_3}_{a_2a_4}\subseteq  P_{a_2a_4}$ and $3 \sqrt[3]{d^*}+d^*=2$.
		Hence, each $5$-edge always gets a color.
		
		By Theorem~\ref{thm-Ramsey}, there exists a subset $S\subseteq [I_1]$ with $|S|=2n+1$ such that all edges $Q$ induced on $S$ have same color. For convenience, we rearrange the indices in $S$ and write $S=[2n+1]$. 
		Suppose that the common color is blue on $S$. Let $I'=\{2,4,\dots, 2n\}\subset S$.
		For any $ijk\ell\in [I']^4$, since the $5$-edge $Q=\{i,j,k, \ell,a\}$ with $a\in S$ and $a>\ell$ is blue, we have $|R^i_{j\ell}\cap T^k_{j\ell}| \ge \frac{\eps}{4}|\mathcal P_{j\ell}|\ge \delta|\mathcal P_{j\ell}|$. 
		Since the induced subhypergraph $\mathcal A_2 \subseteq \mathcal A_1$ with index set $I'$ satisfies $|R^i_{j\ell}\cap T^k_{j\ell}| \ge \delta|\mathcal P_{j\ell}|$ for all $ijk\ell\in [I']^4$, the reverse $\mathcal A'$ of $\mathcal A_2$ on set $I'$ satisfies $|L^{\ell}_{ik}\cap T^{j}_{ik}|\ge \delta|\mathcal P_{ik}|$ for all  $ijk\ell \in [I']^4$, and $\mathcal A'$ is also an induced subhypergraph of $\mathcal A$.
		
		If the common color is red on $S$, we directly choose $I=\{2,4,\dots, 2n\}\subset S$ and the induced subhypergraph $\mathcal A' \subseteq \mathcal A_1$ on the index set $I$. For any $ijk\ell\in [I]^4$, we consider the $5$-edge $Q=\{a,i,j,k, \ell\}$ with $a\in S$ and $a<i$. Since $Q$ is red, we have $|L^{\ell}_{ik}\cap T^{j}_{ik}|\ge \frac{\eps}{4}|\mathcal P_{ik}| \ge \delta|\mathcal P_{ik}|$. 
	\end{proof}
	
	Now we will give an embedding lemma for $(d^*+\eps)$-dense reduced $3$-graphs $\mathcal A$.

	\begin{lemma}\label{embed-two-edge0-lem}
		For any $\eps>0$ and $n\in \mathbb N$, there exists $N\in \mathbb N$ such that the following holds. For every $(d^*+\eps)$-dense $[N]$-reduced $3$-graph $\mathcal A$, there exists an induced subhypergraph $\mathcal A' \subseteq \mathcal A$ with index set $I$ of size $n$ and there exist vertices $x_{ij}, y_{ij}, z_{ij}, y^1_{ij}, z^1_{ij}\in \mathcal P_{ij}$ for all $ij\in [I]^2$ such that $\{x_{ij}, y_{jk}, z_{ik}\}, \{z_{ij}, y^1_{jk}, z^1_{ik}\}\in E(\mathcal A_{ijk})$ for all $ijk\in [I]^3$.
	\end{lemma}
	
	\begin{proof} 
		Given $\eps>0$ and $n\in \mathbb N$, we choose $N$ as follows: First apply Lemma~\ref{left-top-right} with $\eps$ to get $\delta> 0$, then  apply Lemma~\ref{top-lem} with $\delta/2$ and $n$ to get $N_1$, then  apply  Lemma~\ref {right-lem-1} with $\delta/2$ and $N_1$ to get $N_2$,
		then Lemma~\ref{left-lem-1} with $\delta/2$ and $N_2$ to get $N_3$, 
		then Lemma~\ref {right-lem-1}  with $\delta/2$ and $N_3$ to get $N_4$,
		then Lemma~\ref{right-lem} with $\delta$ and $N_4$ to get $N_5$, finally apply Lemma~\ref{left-top-right} with $\eps$, $\delta$ and $N_5$ to get $N$.

		Let $\mathcal A$ be a $(d^*+\eps)$-dense $[N]$-reduced $3$-graph. By Lemma~\ref{left-top-right}, there exists an induced subhypergraph $\mathcal A_5 \subseteq \mathcal A$ with index set $I_5$ of size $N_5$ such that
		\begin{equation}\label{ineq-L-T}
			|L^{\ell}_{ik}\cap T^{j}_{ik}|\ge \delta|\mathcal P_{ik}| \text{~for~all~}  ijk\ell \in [I_5]^4.
		\end{equation}

		By Lemma~\ref{right-lem} applied with $\mathcal A_5$ and sets $L^{\ell}_{ik}\cap T^{j}_{ik}$, there exists an induced subhypergraph $\mathcal A_4\subseteq \mathcal A_5$  with index set $I_4\subseteq I_5$ of size $N_4$, and vertices $z_{ik}\in \mathcal P_{ik}$ for all $ik\in [I_4]^2$ such that 
		\[ 
		z_{ik}\in \bigcap_{i<j<k<\ell,~j\ell \in [I_4]^2} L^{\ell}_{ik}\cap T^{j}_{ik}.
		\]
		For each $ijk\in [I_4]^3$, let 
		\[
		Y_{ijk}=\{y\in \mathcal P_{jk}: \{y, z_{ik}\} \text{~is~contained~in~at~least~} \delta/2 |\mathcal P_{ij}| \text{~hyperedges~of~} \mathcal A_{ijk}\}.
		\] 
		Since $z_{ik}\in T^{j}_{ik}$,  is contained in at least $\delta |\mathcal P_{ij}||\mathcal P_{jk}|$ hyperedges in $\mathcal A_{ijk}$, we have 
		\[
		\delta |\mathcal P_{ij}||\mathcal P_{jk}| \le 
		\delta/2 |\mathcal P_{ij}||\mathcal P_{jk}\setminus Y_{ijk}|+ |Y_{ijk}||\mathcal P_{ij}|,
		\]
		which implies that $|Y_{ijk}|\ge \delta/2|\mathcal P_{jk}|$.
		By Lemma~\ref{right-lem-1} applied with $\mathcal A_4$ and sets $Y_{ijk}$, there exists an induced subhypergraph $\mathcal A_3\subseteq \mathcal A_4$ with index set $I_3\subseteq I_4$ of size $N_3$, and vertices $y_{jk}\in \mathcal P_{jk}$ for all $jk\in [I_3]^2$ such that 
		\[ 
		y_{jk}\in \bigcap_{i<j<k,~i \in I_3} Y_{ijk}.
		\]
		Next, for each $ijk\in [I_3]^3$, let 
		$X_{ijk}=\{x\in \mathcal P_{ij}: \{x,y_{jk}, z_{ik}\}\in E(\mathcal A_{ijk})\}$.
		Clearly, $|X_{ijk}|\ge \delta/2|\mathcal P_{ij}|$.
		By Lemma~\ref{left-lem-1} applied with $\mathcal A_3$ and sets $X_{ijk}$, there exists an induced subhypergraph $\mathcal A_2\subseteq \mathcal A_3$ with index set $I_2\subseteq I_3$ of size $N_2$, and vertices $x_{ij}\in \mathcal P_{ij}$ for all $ij\in [I_2]^2$ such that 
		\[ 
		x_{ij}\in \bigcap_{i<j<k,~k \in I_2} X_{ijk}.
		\]
		
		Until now, we have $x_{ij}, y_{ij}, z_{ij}\in \mathcal P_{ij}$ for all $ij\in [I_2]^2$ and $\{x_{ij}, y_{jk}, z_{ik}\}\in E(\mathcal A_{ijk})$ for all $ijk\in [I_2]^3$.
		Next for each $ijk\in [I_2]^3$, let 
		\[
		Y^1_{ijk}=\{y'\in \mathcal P_{jk}: \{z_{ij}, y'\} \text{~is~contained~in~at~least~} \delta/2 |\mathcal P_{ik}| \text{~hyperedges~of~} \mathcal A_{ijk}\}.
		\] 
		Since $z_{ij}\in L^{k}_{ij}$, is contained in at least $\delta |\mathcal P_{jk}||\mathcal P_{ik}|$ hyperedges in $\mathcal A_{ijk}$, we have 
		\[
		\delta |\mathcal P_{jk}||\mathcal P_{ik}|\le \delta/2 |\mathcal P_{ik}||\mathcal P_{jk}\setminus Y^1_{ijk}|+ |Y^1_{ijk}||\mathcal P_{ik}|,
		\]
		which implies that $|Y^1_{ijk}|\ge \delta/2  |\mathcal P_{jk}|$.
		By Lemma~\ref{right-lem-1} applied with $\mathcal A_2$ and sets $Y^1_{ijk}$, there exists an induced subhypergraph $\mathcal A_1\subseteq \mathcal A_2$ with index set $I_1\subseteq I_2$ of size $N_1$, and vertices $y^1_{jk}\in \mathcal P_{jk}$ for all $jk\in [I_1]^2$ such that 
		\[ 
		y^1_{jk}\in \bigcap_{i<j<k,~i \in I_1} Y^1_{ijk}.
		\]
		Finally for each $ijk\in [I_1]^3$, let 
		$Z^1_{ijk}=\{z'\in \mathcal P_{ik}: \{z_{ij},y^1_{jk}, z^1_{ik}\}\in E(\mathcal A_{ijk}) \}$. 
		Clearly, $|Z^1_{ijk}|\ge \delta/2|\mathcal P_{ik}|$.
		By Lemma~\ref{top-lem} applied with $\mathcal A_1$ and sets $Z^1_{ijk}$, there exists an induced subhypergraph $\mathcal A'\subseteq \mathcal A_1$ with index set $I\subseteq I_1$ of size $n$, and vertices $z^1_{ik}\in \mathcal P_{ik}$ for all $ik\in [I]^2$ such that 
		\[ 
		z^1_{ik}\in \bigcap_{i<j<k,~j \in I} Z^1_{ijk}.
		\]
		Therefore, the $I$-reduced subhypergraph $\mathcal A'$ together with the vertices  $x_{ij}, y_{ij}, z_{ij}, y^1_{ij}, z^1_{ij}\in \mathcal P_{ij}$ for all $ij\in [I]^2$ satisfies
		the conclusion described in the statement of the lemma. 
	\end{proof}

	Finally, let us prove Lemma~\ref{less-than1/4-lem} using Lemma~\ref{embed-two-edge0-lem}.
	\begin{proof}[{\bf Proof of Lemma~\ref{less-than1/4-lem}}]
		Given a $3$-graph $F$ with $f$ vertices satisfying the properties given  in Lemma~\ref{less-than1/4-lem}, we only need to show that for every $\eps>0$, there exists $N\in \mathbb N$ such that every $(d^*+\eps)$-dense $[N]$-reduced $3$-graph can  embed $F$ by Theorem~\ref{thm-turan-reduced}.
		
		We apply Lemma~\ref{embed-two-edge0-lem} with $\eps$ and $f$ to get
		$N$. Let $\mathcal A$ be a $(d^*+\eps)$-dense  $[N]$-reduced $3$-graph. By Lemma~\ref{embed-two-edge0-lem}, there exists an induced subhypergraph $\mathcal A'\subseteq \mathcal A$ 
		with index set $I$ of size $f$ and vertices $x_{ij}, y_{ij}, z_{ij}, y^1_{ij}, z^1_{ij}\in \mathcal P_{ij}$ for all $ij\in [I]^2$ such that $\{x_{ij}, y_{jk}, z_{ik}\}, \{z_{ij}, y^1_{jk}, z^1_{ik}\}\in E(\mathcal A_{ijk})$ for all $ijk\in [I]^3$.

		For brevity, let $I=\{1, \dots,f\}$ and $\phi:V(F)\to I$ satisfying $\phi(v_\ell)=\ell$ for all $\ell \in [f]$. For any pair $\{v_i,v_j\}\in \partial F$ with $1 \le i<j\le f$, we choose the vertex $z_{ij}$ from $\mathcal P_{ij}$ if $\chi (v_i,v_j)$ is \textcolor{blue}{blue}, choose the vertex $y^1_{ij}$ from $\mathcal P_{ij}$ if $\chi (v_i,v_j)$ is \textcolor{violet}{violet}, choose the vertex $z^1_{ij}$ from $\mathcal P_{ij}$ if $\chi (v_i,v_j)$ is \textcolor{red}{red}, choose the vertex $x_{ij}$ from $\mathcal P_{ij}$ if $\chi (v_i,v_j)$ is \textcolor{green}{green}, and choose the vertex $y_{ij}$ from $\mathcal P_{ij}$ if $\chi (v_i,v_j)$ is \textcolor{cyan}{cyan}.
		Since $\{x_{ij}, y_{jk}, z_{ik}\}, \{z_{ij}, y^1_{jk}, z^1_{ik}\}\in E(\mathcal A_{ijk})$ for all $ijk\in [I]^3$,  $\mathcal A'$ can embed $F$.
	\end{proof}

	\section{Applications of Theorem~\ref{main-thm1}}\label{Applications}
	In this section, we will  prove Theorem~\ref{7-k43-} and Theorem~\ref{thm:3-chromatic} using Theorem~\ref{main-thm1}. We start with the following result.
	
	\begin{lemma}\label{low-lemma}
For any even integer $t\ge 4$ and any $3$-graph $G^{a}_t\in \mathcal G^{(3)}_t$, we have $G^{a}_t\in \mathcal F^{\ge 1/4}$.
	\end{lemma}
	
\begin{proof}
Given a 3-graph $G^{a}_t$, assume that there is an ordering $(v_1, \dots, v_t)$ of $V(G^{a}_t)$ and there is a  $2$-coloring $\chi: \partial G^{a}_t \to \{{\color {red}red}, {\color{blue}blue}\}$  such that every $3$-edge $\{v_i, v_j, v_k\}$ of $G^{a}_t$ with $i<j< k$ satisfies
\begin{align}\label{eq:2-coloring}
	\left(\chi(v_i,v_j), \chi(v_j,v_k), \chi(v_i, v_k)\right)=({\color {red}red}, {\color {red}red}, {\color{blue}blue}) \text{~or~} ({\color{blue}blue},  {\color{blue}blue}, {\color {red}red}).
\end{align}

Next, let $a=v_{\ell}$ for some $\ell\in [t]$, and define 
\[
\begin{split}
	& N^{<}_{\color{red}{red}}(v_{\ell})=\{v_i: i<\ell, \chi(v_i,v_{\ell})={\color{red} red}\}, ~~N^{>}_{\color{red}{red}}(v_{\ell})=\{v_i: i>\ell, \chi(v_{\ell},v_i)={\color{red} red}\} \text{~and~}\\
	& N^{<}_{\color{blue}{blue}}(v_{\ell})=\{v_i: i<\ell, \chi(v_i,v_{\ell})={\color{blue} blue}\}, ~~N^{>}_{\color{blue}{blue}}(v_{\ell})=\{v_i: i>\ell, \chi(v_{\ell},v_i)={\color{blue}blue}\}.
\end{split}
\]
By the definition of the above four sets, we get a partition of $V(G^{a}_t)\setminus \{a\}$. Note that since each $3$-edge $e=\{v_{\ell}, v_i, v_j\}\in E(G^{a}_t)$ satisfies the condition~\eqref{eq:2-coloring}, the above four sets are independent sets in the link graph of the vertex $v_{\ell}$. Furthermore, we can check that by the condition~\eqref{eq:2-coloring}:
\[\begin{split}
	& N_{G^{a}_t}(v_{\ell})\cap \{\{x, y\}: x\in N^{<}_{\color{red}{red}}(v_{\ell}), y\in N^{>}_{\color{blue}{blue}}(v_{\ell})\}=\emptyset \text{~and~}\\
	& N_{G^{a}_t}(v_{\ell})\cap \{\{x, y\}: x\in N^{<}_{\color{blue}{blue}}(v_{\ell}), y\in N^{>}_{\color{red}{red}}(v_{\ell})\}=\emptyset.
\end{split}
\]
Therefore, $N^{<}_{\color{red}{red}}(v_{\ell})\cup N^{>}_{\color{blue}{blue}}(v_{\ell})$ and $N^{<}_{\color{blue}{blue}}(v_{\ell})\cup N^{>}_{\color{red}{red}}(v_{\ell})$ are also independent sets in the link graph of the vertex $v_{\ell}$, which implies that $\chi (G)=2$ and it is impossible.
Thus, $G^{a}_t$  does not satisfy the property {\rm ($\clubsuit$)}, i.e., $G^{a}_t\in \mathcal F^{\ge 1/4}$.
	\end{proof}

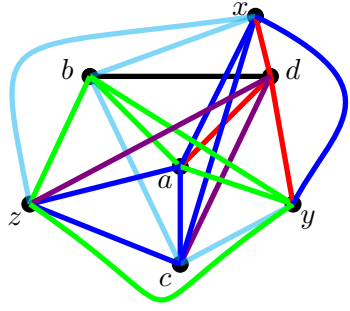
\begin{figure}
		\begin{center}
			\begin{tikzpicture}
				[inner sep=1pt,
				vertex/.style={circle, draw=blue!50, fill=blue!50},
				]
				\filldraw[black] (5,1.5) circle (3pt);
				\node at (4.8,1.3) {$a$};
				\filldraw[black] (5,0.2) circle (3pt);
				\node at (4.8,0) {$c$};
				\filldraw[black] (3.8,2.7) circle (3pt);
				\node at (3.5,2.8) {$b$};
				\filldraw[black] (6.2,2.7) circle (3pt);
				\node at (6.5,2.8) {$d$};
				\draw[line width=2pt,blue,opacity=1]  (5,1.5) --  (5,0.2);
				\draw[line width=2pt,green,opacity=1] (5,1.5) --  (3.8,2.7);
				\draw[line width=2pt,red,opacity=1] (5,1.5) --  (6.2,2.7);
				\draw[line width=2pt,cyan,opacity=0.5] (3.8,2.7) --  (5,0.2);
				\draw[line width=2pt,black,opacity=1] (3.8,2.7) --  (6.2,2.7);
				\draw[line width=2pt,violet,opacity=1] (5,0.2) --  (6.2,2.7);
				\filldraw[black] (6,3.5) circle (3pt);
				\node at (5.8,3.6) {$x$};
				\draw[line width=2pt,blue,opacity=1] (6,3.5) --  (5,1.5);
				\draw[line width=2pt,cyan,opacity=0.5] (6,3.5) --  (3.8,2.7);
				\draw[line width=2pt,blue,opacity=1] (6,3.5) --  (5,0.2);
				\draw[line width=2pt,red,opacity=1] (6,3.5) --  (6.2,2.7);
				\filldraw[black] (6.5,1) circle (3pt);
				\node at (6.7,0.8) {$y$};
				\draw[line width=2pt,green,opacity=1] (6.5,1) --  (5,1.5);
				\draw[line width=2pt,green,opacity=1] (6.5,1) --  (3.8,2.7);
				\draw[line width=2pt,cyan,opacity=0.5] (6.5,1) --  (5,0.2);
				\draw[line width=2pt,red,opacity=1] (6.5,1) --  (6.2,2.7);
				\filldraw[black] (3,1) circle (3pt);
				\node at (2.8,0.8) {$z$};
				\draw[line width=2pt,blue,opacity=1] (3,1) --  (5,1.5);
				\draw[line width=2pt,green,opacity=1] (3,1) --  (3.8,2.7);
				\draw[line width=2pt,blue,opacity=1] (3,1) --  (5,0.2);
				\draw[line width=2pt,violet,opacity=1] (3,1) --  (6.2,2.7);
				\draw[line width=2pt, green, opacity=1] (6.5,1)  .. controls (4,-0.7) and (5.5,-0.7) .. (3,1); 
				\draw[line width=2pt, cyan,opacity=0.5]  (6,3.5)  .. controls (2.5,3) and (2.5,3.5) .. (3,1); 
				\draw[line width=2pt, blue, opacity=1] (6.5,1)  .. controls (7,2) and (8,2.5) .. (6,3.5); 
			\end{tikzpicture}
			\caption{An illustration of  $\partial \hat F^{\star}_7$ with $6$-clorings.}\label{figure-1}
		\end{center}
	\end{figure}

	\begin{proof}[{\bf Proof of Theorem~\ref{7-k43-}}]
		Let $\hat F^{\star}_7$ be the $3$-graph with vertex set $\{a,b,c,d,x,y,z\}$ and edge set $\{abc, abd, acd, xab, xcd, yac, ybd, zad, zbc, xyz\}$.
		Using Theorem~\ref{main-thm1}, we only need to verify that  $
		\hat F^{\star}_7\in \mathcal F^{\ge 1/4}\cap \mathcal F^{\le 1/4}$.  Note first that  $K^{(3)-}_4\sub \hat F^{\star}_7$, we directly obtain that $\hat F^{\star}_7\in \mathcal F^{\ge 1/4}$ by Lemma~\ref{low-lemma}.
		
		Next we will check that $\hat F^{\star}_7$ satisfies the property {\rm ($\spadesuit$)}. 
		Let $\sigma:=(a, y, z, b, x, c, d)$ be an 
		ordering of $V(\hat F^{\star}_7)$.
		Moreover, let  $\chi: \partial \hat F^{\star}_7\to \{\textcolor{red}{red},\textcolor{blue}{blue}, \textcolor{green}{green},\textcolor{violet}{violet},\textcolor{cyan}{cyan}, \textcolor{black}{black}\}$ be a $6$-coloring satisfying the following properties.
		\begin{enumerate}
			\item[{\rm (1)}] All pairs in set $E_1=\{\{a,d\}, \{y,d\}, \{x,d\}\}$ are colored \textcolor{red}{red}; 
			\item[{\rm (2)}] All pairs in set $E_2=\{\{a,z\}, \{a,x\}, \{a,c\}, \{y,x\}, \{z,c\}, \{x,c\}\}$ are colored \textcolor{blue}{blue};
			\item[{\rm (3)}] All pairs in set $E_3=\{\{a,y\},\{a,b\}, \{y,z\}, \{y,b\}, \{z,b\}\}$ are colored \textcolor{green}{green};
			\item[{\rm (4)}] All pairs in set $E_4=\{\{z,d\}, \{c,d\} \}$ are colored \textcolor{violet}{violet}; 
			\item[{\rm (5)}] All pairs in set $E_5=\{ \{y,c\}, \{b,x\}, \{b,c\}, \{z,x\}\}$ are colored \textcolor{cyan}{cyan}; 
			\item[{\rm (6)}] All pairs in set $E_6=\{\{b,d\}\}$ are colored \textcolor{black}{black}. 
		\end{enumerate}
		Then,  under the ordering $\sigma$ and coloring $\chi$,  every $3$-edge of $\hat F^{\star}_7$ satisfies exactly one of the cases given in the property {\rm ($\spadesuit$)}. Indeed, the integer $i^{\star}$ in the property {\rm ($\spadesuit$)} takes $6$ and the hyperedge set $E(\hat F^{\star}_7)$ is partitioned into the following three cases (see Figure~\ref{figure-1}):
		\begin{enumerate}
			\item [(a)] $\left(\chi(i,j), \chi(j,k), \chi(i,k)\right)=(\textcolor{blue}{blue},\textcolor{violet}{violet},\textcolor{red}{red})$ for $ijk\in \{acd, xcd, azd\}$;
			\item [(b)] $\left(\chi(i,j), \chi(j,k), \chi(i,k)\right)=(\textcolor{green}{green},\textcolor{cyan}{cyan},\textcolor{blue}{blue})$ for $ijk\in \{abc, abx, ayc, zbc, yzx\}$;
			\item [(c)] $\left(\chi(i,j), \chi(j,k), \chi(i,k)\right)=(\textcolor{green}{green},\textcolor{black}{black},\textcolor{red}{red})$ for $ijk\in \{abd, ybd\}$.
		\end{enumerate}
	\end{proof}
	
Finally, we will give the proof of Theorem~\ref{thm:3-chromatic} using Lemma~\ref{low-lemma} and Theorem~\ref{main-thm1}. 

\begin{proof}[{\bf Proof of Theorem~\ref{thm:3-chromatic}}]
	Given a 3-graph $G^{a}_t\in \mathcal G^{(3)}_t$ with $t\ge 4$, 
	using Theorem~\ref{main-thm1}, we only need to verify that  $G^{a}_t\in \mathcal F^{\ge 1/4}\cap \mathcal F^{\le 1/4}$. By Lemma~\ref{low-lemma}, we already have $G^{a}_t\in \mathcal F^{\ge 1/4}$. 
	
	Next, we will check that $G^{a}_t\in \mathcal F^{\le 1/4}$.
	Due to $\chi(G)=3$, we consider a partition $\{X_1, X_2, X_3\}$ of $V(G)$ such that $X_i$ is an independent set of $G$ for each $i\in [3]$.
	We consider any ordering $\sigma=(v_1, \dots, v_t)$ of $V(G^{a}_t)$ that satisfies $\{v_1\}=\{a\}$, $\{v_2, \dots, v_{|X_1|+1}\}=X_1$, $\{v_{|X_1|+2}, \dots, v_{|X_1\cup X_2|+1}\}=X_2$, and $\{v_{|X_1\cup X_2|+2}, \dots, v_{t}\}=X_3$.
	Moreover, let  $\chi: \partial G^{a}_t\to \{\textcolor{red}{red},\textcolor{blue}{blue}, \textcolor{green}{green},\textcolor{violet}{violet},\textcolor{cyan}{cyan}, \textcolor{black}{black}\}$ be a $6$-coloring satisfying the following properties:
	\begin{enumerate}
		\item[{\rm (1)}] $\chi (v_1, x)=\textcolor{green}{green}$ for all $\{v_1, x\}\in \partial G^{a}_t$ and $ x\in X_1$;
		\item[{\rm (2)}] $\chi (v_1, x)=\textcolor{blue}{blue}$ for all $\{v_1, x\}\in \partial G^{a}_t$ and $ x\in X_2$;
		\item[{\rm (3)}] $\chi (v_1, x)= \textcolor{red}{red}$ for all $\{v_1, x\}\in \partial G^{a}_t$ and $ x\in X_3$; 
		\item[{\rm (4)}] $\chi (x_1, x_2)=\textcolor{cyan}{cyan}$ for all $\{x_1, x_2\}\in \partial G^{a}_t$ and $(x_1, x_2)\in X_1\times X_2$; 
		\item[{\rm (5)}] $\chi (x_2, x_3)=\textcolor{violet}{violet}$ for all $\{x_2, x_3\}\in \partial G^{a}_t$ and $(x_2, x_3)\in X_2\times X_3$;  
		\item[{\rm (6)}] $\chi (x_1, x_3)=\textcolor{black}{black}$  for all $\{x_1, x_3\}\in \partial G^{a}_t$ and $(x_1, x_3)\in X_1\times X_3$.
	\end{enumerate}
	Finally, the integer $i^{\star}$ in the property {\rm ($\spadesuit$)}  takes $|X_1|+|X_2|+1$.
	Clearly, under the ordering $\sigma$ and coloring $\chi$,  every $3$-edge of $G^{a}_t$ satisfies exactly one of the cases given in the property {\rm ($\spadesuit$)}. 
	Therefore, $G^{a}_t$ satisfies the property {\rm ($\spadesuit$)}, i.e., $G^{a}_t\in \mathcal F^{\ge 1/4}$.
\end{proof}

\noindent {\bf Foundations}
Guanghui Wang's research is supported by the Natural Science Foundation of China  (No. 12231018).
Wenling Zhou's research is supported by  the Natural Science Foundation of China (No. 12401457),  the China Postdoctoral Science Foundation (No. 2024M761780)  and Natural Science Foundation of Shandong Province (No. ZR2024QA067).

\smallskip

\noindent {\bf Author Contributions} Hao Li, Hao Lin, Guanghui Wang and Wenling Zhou finished the theoretical proof together. All authors have read and approved the final manuscript.
\smallskip

\noindent {\bf Conflicts of interest} The authors declare no conflict of interest. The funding sponsors had no role in the design of the study; in the collection, analysis or interpretation of data; in the writing of the manuscript; nor in the decision to publish the results.

	\bibliographystyle{abbrv}
	\bibliography{revised-ref}
\end{document}